\documentclass[12pt]{article}

\usepackage{graphicx}
\usepackage{mathptmx}
\usepackage{latexsym,amsmath,amssymb,amsfonts,amsthm}
\newcommand{\R}{\mathbb{R}}

\newcommand{\sm}[1]{\mbox{\small $#1$}}

\newcommand{\La}[1]{\mbox{\Large $#1$}}

\newcommand{\oM}{\bar{M}}

\newtheorem{theorem}{Theorem}[section]
\newtheorem{lemma}{Lemma}[section]
\newtheorem{corollary}{Corollary}[section]
\newtheorem{proposition}{Proposition}[section]

\begin{document}
\setcounter{page}{1}
\title{Mean curvature flow of spacelike graphs}
\author{Guanghan Li$^{1,2\dag}$ and  Isabel M.C.\ Salavessa$^{2\ddag}$}
\pagestyle{myheadings}
\markright{\sl\hfill Li and Salavessa,  Mean curvature flow of spacelike
graphs \hfill}
\date{}
\protect\footnotetext{\!\!\!\!\!\!\!\!\!\!\!\!\! {\bf MSC 2000:}
Primary:53C21, 53C40 Secondary: 58D25, 35K55.
$\dag$ Partially supported by NSFC (No.10501011) and by
Funda\c{c}\~{a}o Ci\^{e}ncia e Tecnologia (FCT) through a FCT
fellowship SFRH/BPD/26554/2006. \\
$\ddag$ Partially supported by FCT through the Plurianual of CFIF
and POCI-PPCDT/MAT/60671/2004.
}
\maketitle ~~~\\[-15mm]
\begin{quotation}\noindent
\baselineskip .3cm{
{\footnotesize $^1$ School of Mathematics and Computer
Science, Hubei University, Wuhan, 430062, P. R. China,
$\dag$ e-mail: liguanghan@163.com}\\[2mm]
{
\footnotesize $^2$ Centro de F\'{\i}sica das Interac\c{c}\~{o}es
Fundamentais, Instituto Superior T\'{e}cnico, Technical University
of Lisbon, Edif\'{\i}cio Ci\^{e}ncia, Piso 3, Av.\ Rovisco Pais,
1049-001 Lisboa, Portugal;~~
$\ddag$ e-mail: isabel.salavessa@ist.utl.pt}}
\end{quotation}
\begin{abstract} \noindent 
We prove the mean curvature flow of
a spacelike graph in $(\Sigma_1\times \Sigma_2, g_1-g_2)$ of a map
 $f:\Sigma_1\to \Sigma_2$ from a closed  Riemannian
manifold $(\Sigma_1,g_1)$ with  $Ricci_1> 0$ 
to a complete
Riemannian manifold
$(\Sigma_2,g_2)$ with bounded curvature tensor and derivatives,
and with sectional curvatures satisfying $K_2\leq K_1$,
remains a spacelike graph, exists for all time, and converges to a slice
at infinity. 
We also  show, with no need of the assumption
$K_2\leq K_1$, that if $K_1>0$, or if $Ricci_1>0$ and 
$K_2\leq -c$, $c>0$  constant, 
any map $f:\Sigma_1\to \Sigma_2$
is trivially homotopic provided $f^*g_2<\rho g_1$ where
$\rho=\min_{\Sigma_1}K_1/\sup_{\Sigma_2}K_2^+\geq 0$, in case $K_1>0$,
and $\rho=+\infty$ in case $K_2\leq 0$. This largely extends some known
results for $K_i$ constant and $\Sigma_2$ compact,
 obtained using the Riemannian structure of
$\Sigma_1\times \Sigma_2$, and also
 shows how regularity theory on the mean curvature flow is
simpler and more natural in pseudo-Riemannian setting then
in the Riemannian one.
\end{abstract}

\section{Introduction}
\label{Introduction}
\renewcommand{\thesection}{\arabic{section}}
\renewcommand{\theequation}{\thesection.\arabic{equation}}
\setcounter{equation}{0}
\setcounter{theorem}{0}
\setcounter{corollary}{0}
\setcounter{proposition}{0}
\setcounter{lemma}{0}

Let $M$ be a smooth manifold of dimension $m\ge 2$, and $F_0:
M\rightarrow \oM$  a smooth submanifold immersed into a
$(m\!+\!n)$-dimensional Riemannian or pseudo-Riemannian manifold $(\oM,\bar{g})
$. The mean curvature flow is a smooth family of immersions
$F_t=F(\cdot, t):M\rightarrow \oM$ evolving according to
\begin{eqnarray}\left\{\begin{array}{lll}
\frac d{dt}F(x, t)&=&{H}(x,t), \quad x\in M,\\
[2mm]   F(\cdot, 0)&=&F_0, \end{array}\right.
\end{eqnarray}
where $H$ is the mean curvature vector of $M_t=(M,g_t=F^*_t\bar{g})=F_t(M)$.

The mean curvature flow of hypersurfaces (i.e. (1.1) with $n=1$)
in a Riemannian manifold has been extensively studied in the last
two decades. Recently, mean curvature flow of submanifolds with
higher co-dimensions has been paid more attention,
 see \cite{cl,s,w1,w2,x2} for
example. In \cite{w2}, the graph mean curvature
flow is studied in Riemannian product manifolds, and it is proved long-time
existence and convergence of the  flow under
suitable conditions.

When $\oM$ is a pseudo-Riemannian manifold, (1.1) is  the
mean curvature flow of spacelike submanifolds. This flow
for spacelike hypersurfaces has also been  strongly studied, see
\cite{aa,e,eh,g} and references therein. To our knowledge, very little
is known on mean curvature flow in higher codimensions except
in a flat space  $\mathbb{R}^{n+m}_n$ \cite{x2}.
In this paper,
 we partially
 apply  Wang's approach \cite{w2} of using the mean curvature flow
in a Riemannian product space
to deform a map between Riemannian manifolds, but we use the
pseudo-Riemannian structure of the product.
As a result we obtain a reformulated and largely extended version
of the the main results in \cite{w2} and most of \cite{tw},
to  non constant sectional curvatures $K_i$,
and applied to a set of maps satisfying a less restrictive condition,
after using a simple argument of rescaling the metric
in the target space
$\Sigma_2$,
in  a convenient way.
 The pseudo-Riemannian structure turns out
to give a simpler and more natural tool to provide an existence result
on the deformation of a map to a
totally geodesic one or to a constant one
by some curvature flow under quite weak curvature conditions.
In
\cite{w2} it is necessary to use  White's regularity theorem
\cite{wh}, where a monotonicity formula due to Huisken
\cite{h3} plays a fundamental role,
 to detect possible singularities of the mean curvature flow, 
while in pseudo-Riemannian
case, because of good signature in the evolution equations, we have
better regularity, and therefore we require fewer restrictions on
the curvatures of $\Sigma_1$ and $\Sigma_2$ in the main theorem
\ref{theorem 1.1} as well on the map $f$ itself in theorem
\ref{theorem 1.2}. We also believe that this ``pseudo-Riemannian" trick
may be applied to some other geometric evolution equations to obtain
the convergence of the flow in a more efficient way.

Let $\oM=\Sigma_1\times \Sigma _2$ be a product manifold of two
Riemannian manifolds $(\Sigma_1, g_1)$ of dimension $m$
and $(\Sigma_2, g_2)$ of dimension $n$, with
the pseudo-Riemannian metric $\bar{g}=g_1-g_2$, where $\Sigma
_i, i=1,2$ has sectional curvature
$K_i $ and Ricci tensor $Ricci_i$. Assume $M$ is a spacelike graph
$$M=\Gamma_f=\{(p,f(p)):
p\in \Sigma_1\}$$
  of a smooth map $f:\Sigma_1\rightarrow\Sigma
_2$, with induced metric $g$. The graph map,
$\Gamma_f:\Sigma_1\rightarrow M$, $\Gamma_f(p)=(p,f(p))$,
identifies isometrically $(M,g)$
with $\Sigma_1$ with the graph metric $g_1-f^*g_2$.
$M$ is a slice if $f$ is a constant map.
The  hyperbolic angle  $\theta$ can be defined by
(see \cite{alal}, \cite{ls})
$$\cosh\theta=\frac{1}{\sqrt{det(g_1-f^*g_2)}},$$
where the determinant is defined with respect to the metric $g_1$.
The angle $\theta$ measures the deviation from a spacelike
submanifold to a slice. If this angle
 is bounded the metrics $g_1$  and $g=g_1-f^*g_2$  of
$\Sigma_1$ are equivalent. In this case,
  $(\Sigma_1,g_1)$ is compact iff  $(M,g)$ is so.
 The
following is the main theorem in this paper.

\begin{theorem} \label{theorem 1.1} Let $f$ be a smooth map from $\Sigma _1$
to $\Sigma _2$ such that $F_0: M\rightarrow \oM$ is a compact
spacelike graph of $f$.  We assume $(\Sigma_1,g_1)$ closed
of dimension $m\geq 2$, $(\Sigma_2,
g_2)$ complete of dimension $n\geq 1$,  $Ricci_1(p)\geq 0$ and 
$K_1(p)\ge K_2(q)$ for any $p\in
\Sigma _1$, $q\in \Sigma _2$ and  the curvature
tensor $R_2$ of $\Sigma_2$ and all its covariant derivatives are bounded.
 Then:\\[3mm]
$(1)$~ The mean curvature flow $(1.1)$ of the spacelike
graph of $f$ remains a spacelike graph of a map
$f_t:\Sigma_1\rightarrow \Sigma_2$ and exists for all time.\\[2mm]
$(2)$~ If  $\Sigma_2$ is also compact there is a sequence $t_n\to +\infty$
such that the  flow converges at infinity to a spacelike graph of a totally
geodesic map, and if $Ricci_1(p)>0$ at some point $p\in \Sigma_1$, the sequence
converges to a slice.\\[2mm]
$(3)$~If $Ricci_1>0$ everywhere, all the flow converges  to a unique slice.
\end{theorem}
\noindent
We observe that in (3) we do not need the compactness assumption of $\Sigma_2$.
We also note that the condition $Ricci_1\geq 0$ and $K_1\geq K_2$ means that
at a point $p\in \Sigma_1$, if $K_1(P)<0$ at some two-plane $P$ of
$T_p\Sigma_1$, then we have to require $K_2<0$ everywhere.
In case (3) the flow defines a homotopy $f_t(\phi_t(p))$ from
the initial map $f_0=f$ to the final
constant map $f_{\infty}$,
where $\phi_t=\pi_1\circ F_t$ is a smooth diffeomorphic endomorphism 
of $\Sigma_1$, that at $t=0$ gives the identity map.
We shall prove that the deformation process
is also valid without assuming $K_1\ge K_2$.
\begin{theorem} \label{theorem 1.2} Suppose $(\Sigma _1, g_1)$ and
$(\Sigma _2, g_2)$ are two complete Riemannian manifolds 
of dimensions $m\geq 2$ and $n\geq 1$ respectively,
$\Sigma_1$ closed,  $K_1>0$
everywhere, or  $Ricci_1>0$ and $K_2\leq -c$ with $c>0$ constant,
and the curvature tensor of $\Sigma_2$
and all its covariant derivatives are bounded.
Then there exists a constant $\rho\geq 0$, depending only
on $\min\, K_1$ and on $\sup\, K_2^+$, such that any smooth map
$f: \Sigma _1 \to \Sigma _2$ satisfying
$f^*g_2<\rho g_1$  can be homotopically
 deformed into a constant
map.
\end{theorem}
\noindent
The constant $\rho\in (0,+\infty]$ can be taken equal to
$ {\min_{\Sigma_1}K_1}/{\sup_{\Sigma_2}{K_2^+}}$,
where $K_2^+=\sup\{0,K_2\}$, in case $K_1>0$, and 
 equal to $+\infty $ in case $Ricci_1>0$ and  $K_2\leq -c$.
Recall that by the Cartan-Hadmard theorem, if  $K_2\leq 0$,
the universal cover of $\Sigma_2$ is diffeomorphic to a
Euclidean space, and in particular $\pi_m(\Sigma_2)=0$ for all $m>0$.
If $\Sigma_1$ is the $m$-sphere, and  $K_2\leq 0$ then $\rho=+\infty$ and  the
previous corollary gives a new proof of this classical result.
\begin{corollary}  \label{corollary 1.1}
If  $\Sigma_1$ is closed with $K_1>0$ everywhere,
$K_2\leq 0$ and for all $k\geq 0$, $\nabla^k R_2$ is bounded, 
 then any map $f:\Sigma_1\to \Sigma_2$
is homotopically trivial.
\end{corollary}
We also may apply the previous theorem to obtain a reformulated
version of the main result in \cite{w2}:

\begin{corollary} \label{corollary 1.2} (\cite{w2})
If $\, \Sigma_1$ and $\Sigma_2$ are
compact with constant sectional curvatures $\tau_1$ and $\tau_2$ satisfying
$\tau_1\geq |\tau_2|$ , $\tau_1+\tau_2>0$, and  if $det(g_1+f^*g_2)<2$, then
$\Gamma_f$ can be deformed by a family of graphs to the one of a constant
map.
\end{corollary}
\noindent
The condition $det(g_1+f^*g_2)<2$ implies $\Gamma_f$ is a spacelike
submanifold
for the pseudo-Riemannian structure of $\Sigma_1\times \Sigma_2$.
The converse may not hold, so spacelike graph is a less restrictive
condition.
In \cite{tw},
corollary \ref{corollary 1.2}
 is generalized, under the same constant curvature conditions,
 to the case of area decreasing maps
that is a slightly less restrictive condition than of a spacelike
graph   $f^*g_2<g_1$. For such maps the eigenvalues $\lambda_i^2$ of
$f^*g_2$ satisfy $\lambda_i\lambda_j<1$ for $i\neq j$. Thus,
$f^*g_2$ may have one and only one eigenvalue (counting with
multiplicity) greater than or equal to one. If $n\geq 2$ area
decreasing maps are spacelike iff the largest eigenvalue $\lambda_1$
is also smaller than one. In this case we also recover the main
theorem of \cite{tw}. If $\Sigma_2$ is one-dimensional, any map $f$
satisfies such condition, and the result can be obtained from
theorem \ref{theorem 1.2}, since $K_2=0$ holds in this case. This is a
particular case of $K_2\leq 0$ stated above.\\[3mm]
We consider the $\phi$-energy functional acting
on smooth maps $f:(\Sigma_1,g_1)\to (\Sigma_2,g_2)$, 
$$E_{\phi}(f)
=\int_{\Sigma_1}\phi(\lambda_1^2,\ldots ,\lambda_m^2)d\mu_1,$$
where $d\mu_1$ means here the volume element of $(\Sigma_1, g_1)$,
$\phi$ is a symmetric nonnegative continuous function on the eigenvalues  
$\lambda=\lambda_i^2$ 
of $f^*g_2$ satisfying $\phi(\lambda)=0$ if and only if $\lambda=0$
and $\phi(\lambda)\leq C\|\lambda\|^{\tau}$, for some constants $C,\tau>0$.
When $\phi(\lambda)=\lambda_1^2+\ldots +\lambda_m^2$, 
we have the usual energy functional
whose critical points are the harmonic maps.
As a corollary of theorem 1.2 we obtain:
\begin{corollary} Under the same curvature  conditions of $(\Sigma_i,g_i)$
given in theorem 1.2, if $f:(\Sigma_1,g_1)\to (\Sigma_2, g_2)$
 minimizes the $\phi$-energy functional in its homotopy class, and
if $f^*g_2\leq  \rho g_1$, 
then $f$ is constant. 
\end{corollary}
The rest of this paper is organized as follows. In section 2, we
derive the elementary formulae for the geometry of spacelike
submanifolds in a pseudo-Riemannian manifold. Section 3 is devoted
to spacelike submanifolds in pseudo-Riemannian product manifolds in
our setting. Evolution equations of different geometric quantities
are given in section 4.  In section 5  we 
 prove part of theorem 1.1(1), and
in section 6  we obtain long-time existence using elliptic
Schauder theory
and prove the existence of a convergent sequence of the flow. The
 use of the Bernstein-type results obtained in \cite{alal,ls} leads to
theorem 1.1(2).
In section 7 we consider the particular case
 $Ricci_1>0$ everywhere, and $\Sigma_2$ not necessarily compact, 
 and prove 
the convergence of all the flow. In this section we also prove
theorem \ref{theorem 1.2} and corollaries 1.2 and 1.3.

\section{Geometry of spacelike submanifolds}
\renewcommand{\thesection}{\arabic{section}}
\renewcommand{\theequation}{\thesection.\arabic{equation}}
\setcounter{equation}{0}
\setcounter{theorem}{0}
\setcounter{corollary}{0}
\setcounter{proposition}{0}
\setcounter{lemma}{0}

 Let $\oM$ be an
$(m+n)$-dimensional pseudo-Riemannian manifold, and $\bar{g}$ the
non-degenerate metric on $\oM$, which is of index $n$. Denote
by $\bar{\nabla}$ the connection on $\oM$, and we convention
that the  curvature tensor $\bar{R}$ is defined by
$\bar{R}(X, Y)Z=\bar{\nabla}_X\bar{\nabla}_YZ
-\bar{\nabla}_Y\bar{\nabla}_XZ-\bar{\nabla}_{[X,Y]}Z,$
and
$\bar{R}(X, Y,Z, W)=\bar{g}(\bar{R}(Z, W)Y, X),$
for any  tangent vector fields $X, Y, Z$ and $W$ of
$\oM$.
Suppose $F: M\rightarrow \oM$ is a $m$-dimensional spacelike
submanifold immersed into $\oM$, i.e. the induced metric of
$M$ is positive definite. For any tangent vector fields $X, Y$ of
$M$ and $V$  a time-like normal vector,
$$\bar {\nabla }_XY=\nabla _XY+B(X, Y),~~~~~~
\bar {\nabla }_XV=\nabla ^{\bot}_XV-A_VX,$$
where  $\nabla $ is the induced
connection on $M$, and $\nabla ^{\bot}_XV=(\bar{\nabla}_XV)^{\bot}$
the normal connection in the normal bundle $NM$, and $B$ and $A$
are the second fundamental form and the Weingarten transformation,
respectively,
$g(A_V(X),Y)=\bar{g}(V,B(X,Y))$.

We choose  orthonormal frame fields $\{e_1, \cdots, e_{m+n}\}$ of
$\oM$, such that when restricting to $M$, $\{e_1, \cdots, e_m\}$
is a tangent frame field, and $\{e_{m+1}, \cdots, e_{m+n}\}$ is
a normal frame field. We make use of the indices range, $1\le i,
j, k, \cdots, \le m$, $m+1\le \alpha, \beta, \cdots, \le m+n$, and
$1\le a, b, c, \cdots, \le m+n$. Let $\theta ^1, \cdots, \theta
^{m+n}$ be the dual frame fields of $\{e_a\}$. Then the structure
equations of $\oM$ are given by
\begin{eqnarray*}
d\theta ^a=-\sum_b \theta _b^a\wedge \theta ^b,~~~~~~~~
d\theta _b^a=-\sum _c \theta _c^a\wedge \theta _b^c+\Phi _b^a,
\end{eqnarray*}
where $\Phi _b^a=\frac
12 \sum _{c,d}\bar{R}^a_{bcd}\theta ^c\wedge \theta ^d$
are the curvature forms, and $\theta
^b_a$ the connection forms satisfying
$\sum_c\bar{g}_{ac}\theta _b^c+\bar{g}_{cb}\theta _a^c=d\bar{g}_{ab}=0$.
Let $\omega ^a =F^*\theta ^a$, $\omega _a^b=F^*\theta _a^b$. Then
restricting to $M$, we have
$\omega ^{\alpha }=0$, and $\sum _{i}\omega _i^{\alpha }\wedge \omega ^i=0$.
By Cartan lemma,
\begin{equation} \label{eq 2.1}
\omega _i^{\alpha }=\sum
_{j}h_{ij}^{\alpha}\omega ^j, \quad h_{ij}^{\alpha
}=h_{ji}^{\alpha},
\end{equation}
where $h_{ij}^{\alpha}$ are the components of the second fundamental
form, that is $B(e_i,
e_j)=\sum_{\alpha}h_{ij}^{\alpha}e_{\alpha}$.
Since the normal vectors are time-like,
 the following relations hold
$$\bar{g}(B(e_i, e_j), e_{\alpha})= \bar{g}(A_{\alpha }e_i,
e_j)=-h_{ij}^{\alpha}.$$
The structure equations of $M$ are then given by
$$\begin{array}{ccl}
d\omega ^i &=&-\sum_j \omega _j^i\wedge \omega ^j,\\
d\omega  _j^i&=&-\sum _k \omega _k^i\wedge \omega _j^k+\Omega _j^i,\\
d\omega _{\beta}^{\alpha} &=&-\sum _{\gamma}\omega _{\gamma}^{\alpha}
\wedge \omega _{\beta}^{\gamma}+\Omega _{\beta }^{\alpha},
\end{array}$$
where $\nabla_{e_j}e_i=\sum_k \omega^k_i(e_j)e_k$,
$\nabla^{\bot}_{e_j}e_{\alpha}=\sum_{\alpha}\omega^{\beta}_{\alpha}(e_j)
e_{\beta}$, with
$\omega^{\beta}_{\alpha}+\omega^{\alpha}_{\beta}=\omega^k_i+\omega^i_k=0$,
and
$\Omega ^i_j=\frac 12\sum _{k,l}R_{jkl}^i\omega ^k\wedge \omega ^l$,
$\sum_i R^i_{jkl}e_i=
R(e_k,e_l)e_j$ is the curvature form of $M$,
and $\Omega _{\beta }^{\alpha}=\frac 12\sum _{k,l}R_{\beta
kl}^{\alpha}\omega ^k\wedge \omega ^l$,
$R^{\bot}(e_j,e_k)e_\beta=\sum_{\alpha}R^{\alpha}_{\beta j k}e_{\alpha}$
is the normal curvature form.
 Setting
$\bar{R}(e_c,e_d)e_b=\sum_a \bar{R}^a_{bcd}e_a$
 we have the  Gauss equation
$$R_{jkl}^i=\bar{R}_{jkl}^i-\sum _{\alpha}(h_{ik}^{\alpha}
h_{jl}^{\alpha}-h_{il}^{\alpha}h_{jk}^{\alpha}),$$
and the normal curvature of $M$ is
 given by the Ricci equation
$$R_{\beta kl}^{\alpha}=\bar{R}_{\beta kl}^{\alpha}-\sum_i(h_{ki}^{\alpha}h_{li}^{\beta}-h_{li}^{\alpha}h_{ki}^{\beta}).$$
The tensor given by $\nabla_{Z}B\, (X,Y)=\nabla^{\bot}_Z(B(X,Y))
-B(\nabla_{Z}X,Y)-B(X,\nabla_{Z}Y)$ is
the covariant derivative of $B$.
The components of $\nabla_{e_k}B (e_i,e_j)=
\sum_{\alpha}h^{\alpha}_{ij,k}e_{\alpha}$ satisfy
\begin{equation}\label{eq 2.2}
\sum _k h_{ij,k}^{\alpha}\omega ^k=dh_{ij}^{\alpha}-\sum
_kh_{kj}^{\alpha}\omega _i^k-\sum _kh_{ik}^{\alpha}\omega _j^k+\sum
_{\beta}h_{ij}^{\beta}\omega _{\beta}^{\alpha}.
\end{equation}
Consider $\tilde{R}:\wedge^2TM\rightarrow L(TM;NM)$ the restriction
of $\bar{R}$, defined for $X,Y,Z\in T_pM$ and $U\in NM_p$,
$\bar{g}(\tilde{R}(X,Y)Z,U)=\bar{g}(\bar{R}(X,Y)Z,U).$
Then the components of $\tilde{R}$ are just
$\bar{R}_{ijk}^{\alpha}=\bar{R}_{aijk}\bar{g}^{a\alpha}=-\bar{R}_{\alpha
ijk}$.
Differentiating both sides of (\ref{eq 2.1}) and applying the structure equations
we have the Codazzi equation
$$h_{ij,k}^{\alpha}-h_{ik,j}^{\alpha}=
-\left((\bar{R}(e_j,e_k)e_i)^{\bot}\right)^{\alpha}=
-\bar{R}_{ijk}^{\alpha}.$$
The mean curvature of $F$ is denoted by
$H =trace\, B= \sum_{\alpha}H^{\alpha}e_{\alpha}$,
$H^{\alpha}=\sum_{i}h^{\alpha}_{ii}$.
The tensor defined by $\nabla^{2}_{X,Y}B \, (Z,W)=(\nabla_Y(\nabla_XB)
-\nabla_{\nabla_{Y}X}B)(Z,W)$ is the
  second covariant derivative of $B$.
The components of $\nabla^{2}_{e_k,e_l}B \, (e_i,e_j)= \sum_{\alpha}
h_{ij,kl}^{\alpha}e_{\alpha}$,  satisfy
\begin{eqnarray*}\sum _l h_{ij,kl}^{\alpha}\omega ^l=dh_{ij,k}^{\alpha}-\sum
_lh_{lj,k}^{\alpha}\omega _i^l-\sum _lh_{il,k}^{\alpha}\omega
_j^l-\sum _lh_{ij,l}^{\alpha}\omega _k^l+\sum
_{\beta}h_{ij,k}^{\beta}\omega _{\beta}^{\alpha}.\nonumber
\end{eqnarray*}
Differentiation of (\ref{eq 2.2}) and use of the structure equations we have
\begin{eqnarray}\label{eq 2.3}
h_{ij,kl}^{\alpha}-h_{ij,lk}^{\alpha}=\sum
_rh_{ir}^{\alpha}R_{jkl}^r+\sum _rh_{rj}^{\alpha}R_{ikl}^r-\sum
_{\beta}h_{ij}^{\beta}R_{\beta kl}^{\alpha}.
\end{eqnarray}
In order to compute the Laplacian of the second fundamental form, we
have to relate the covariant derivatives
$(\bar{\nabla}_{e_s} \bar{R}\, (e_j,e_k)e_i)^{\alpha}$,
with
$\nabla_{e_s}\tilde{R}\, (e_j,e_k)e_i=
\sum_{\alpha}\tilde{R}^{\alpha}_{ijk,s}e_{\alpha}$,
where $\nabla\tilde{R}$ is the covariant derivative considering
the connection of the normal bundle. We have
\begin{eqnarray*}(\bar{\nabla} _l\bar{R})_{ijk}^{\alpha}=
\tilde{ R}_{ijk,l}^{\alpha}-\sum_ {\beta}\bar{R}_{\beta jk}^{\alpha}
h_{il}^{\beta}
-\sum _{\beta}\bar{R}_{i\beta k}^{\alpha}h _{jl}^{\beta}-\sum
_{\beta}\bar{R}_{ij\beta }^{\alpha}h _{kl}^{\beta}+\sum
_{r}\bar{R}_{ijk}^{r}h _{rl}^{\alpha}.\nonumber
\end{eqnarray*}
Using Codazzi's equation (4 times), we obtain
\begin{eqnarray*}
h_{ij,ks}^{\alpha}
&=& h_{ik,js}^{\alpha}-\tilde{R}^{\alpha}_{ijk, s}
\end{eqnarray*}
and so, using this equation again and the commutation formula
(\ref{eq 2.3}), we get
\begin{eqnarray*}
h_{ij,kk}^{\alpha} &=&h_{ki,jk}^{\alpha}-\tilde{R}^{\alpha}_{ijk, k}\\
&=&h_{ki,kj}^{\alpha}+ \sum_r h^{\alpha}_{kr}R^r_{ikj}
+\sum_r h^{\alpha}_{ri}R^r_{kkj}-
\sum_{\beta}h^{\beta}_{ki}R^{\alpha}_{\beta jk}-\tilde{R}^{\alpha}_{ijk, k}\\
&=&h_{kk,ij}^{\alpha}-\tilde{R}^{\alpha}_{kik,j} +\sum_r
h^{\alpha}_{kr}R^r_{ikj} +\sum_r h^{\alpha}_{ri}R^r_{kkj}-
\sum_{\beta}h^{\beta}_{ki}R^{\alpha}_{\beta
jk}-\tilde{R}^{\alpha}_{ijk, k}.
\end{eqnarray*}
Thus,
\begin{eqnarray*}
{\sum}_kh_{ij,kk}^{\alpha}&=&{\sum}
_k(h_{ik,j}^{\alpha}-\tilde{R}_{ijk}^{\alpha})_{,k}\\
&=& {\sum}_k (h_{ki,jk}^{\alpha}-\tilde{R}_{ijk,k}^{\alpha})\nonumber\\
&=&{\sum}_k(~h_{ki,kj}^{\alpha}+\sm{\sum}_r h_{kr}^{\alpha}R_{ijk}^r
+\sm{\sum}_rh_{ri}^{\alpha}R_{kjk}^r-\sm{\sum}_{\beta}h_{ki}^{\beta}R_{\beta
jk}^{\alpha}-\tilde{R}_{ijk,k}^{\alpha})\nonumber\\
&=&{\sum}_k(
h_{kk,ij}^{\alpha}-\tilde{R}_{kik,j}^{\alpha}
+\sm{\sum}_r h_{kr}^{\alpha}R_{ijk}^r
+\sm{\sum}_r h_{ri}^{\alpha}R_{kjk}^r-\sm{\sum}_{\beta}
h_{ki}^{\beta}R_{\beta
jk}^{\alpha}-\tilde{R}_{ijk,k}^{\alpha}).
\end{eqnarray*}
The Laplacian of $B$ is the symmetric $NM$-valued 2-tensor of $M$,
$\Delta B=trace \nabla^{2}_{\cdot,\cdot}B$, that is
$(\Delta B(e_i,e_j))^{\alpha}=\sm{\sum}_k h^{\alpha}_{ij,kk}=
"\Delta  h^{\alpha}_{ij}".$
Then we have
\begin{eqnarray*}
\lefteqn{(\Delta B(e_i,e_j))^{\alpha}=\Delta h_{ij}^{\alpha}=}\\
 &=&
H_{,ij}^{\alpha}+ \sm{\sum}_k(
-\tilde{R}_{kik,j}^{\alpha}-\tilde{R}_{ijk,k}^{\alpha}
+\sm{\sum}_rh_{kr}^{\alpha}R_{ijk}^r+
\sm{\sum}_rh_{ri}^{\alpha}R_{kjk}^r-\sm{\sum}_{\beta}h_{ki}^{\beta}R_{\beta
jk}^{\alpha})\nonumber\\
&=&H_{,ij}^{\alpha}+ \sm{\sum}_k\left(-(\bar{\nabla}_j
\bar{R})_{kik}^{\alpha}-\sm{\sum}_{\beta}\bar{R}_{\beta
ik}^{\alpha}h_{kj}^{\beta}-\sm{\sum}_{\beta}\bar{R}_{k\beta
 k}^{\alpha}h_{ij}^{\beta}-\sm{\sum}_{\beta}\bar{R}_{ki\beta
}^{\alpha}h_{kj}^{\beta}+\sm{\sum}_{l}\bar{R}_{kik}^{l}
h_{lj}^{\alpha}\right.\\
&&~~~-(\bar{\nabla}_k\bar{R})_{ijk}^{\alpha}-\sm{\sum}_{\beta}\bar{R}_{\beta
jk}^{\alpha}h_{ik}^{\beta}-\sm{\sum}_{\beta}\bar{R}_{i\beta
 k}^{\alpha}h_{jk}^{\beta}-\sm{\sum}_{\beta}\bar{R}_{ij\beta
}^{\alpha}h_{kk}^{\beta}
+\sm{\sum}_{l}\bar{R}_{ijk}^{l}h_{lk}^{\alpha}\nonumber\\
&&~~~+\sm{\sum}_{r}h_{kr}^{\alpha}[\bar{R}_{ijk}^r-
\sm{\sum}_{\beta}(h_{rj}^{\beta}h_{ik}^{\beta}-h_{rk}^{\beta}h_{ij}^{\beta})]
+\sm{\sum}_{r}h_{ri}^{\alpha}[\bar{R}_{kjk}^r-
\sm{\sum}_{\beta}(h_{rj}^{\beta}h_{kk}^{\beta}
-h_{rk}^{\beta}h_{kj}^{\beta})]\nonumber\\
&&\left.~~~-\sm{\sum}_{\beta}h_{ki}^{\beta}[\bar{R}_{\beta
jk}^{\alpha}-\sm{\sum}_{l}(h_{jl}^{\alpha}h_{kl}^{\beta}
-h_{kl}^{\alpha}h_{jl}^{\beta})]\right).
\end{eqnarray*}
Using  the first Jacobi identity
$\bar{R}^{\alpha}_{ij\beta}=-\bar{R}^{\alpha}_{\beta ij}
-\bar{R}^{\alpha}_{j\beta i}$,
and that $\sum_{ij}h^{\alpha}_{ij}\bar{R}^{\alpha}_{\beta ij}=0$,
we have
\begin{eqnarray*}
-\langle B,\Delta B\rangle =
\sm{\sum}_{ij\alpha}h_{ij}^{\alpha}\Delta h_{ij}^{\alpha}&=&
\sm{\sum}_{ij\alpha}\La{(}~
h_{ij}^{\alpha}H_{,ij}^{\alpha}-h_{ij}^{\alpha}[\sm{\sum}_k(\bar{\nabla}
_j\bar{R})_{kik}^{\alpha}+\sm{\sum}_k(\bar{\nabla}
_k\bar{R})_{ijk}^{\alpha}]\nonumber\\
&&~~~+\sm{\sum}_{k\beta} (4\bar{R}_{\beta
ki}^{\alpha}h_{kj}^{\beta}h_{ij}^{\alpha}-\bar{R}_{k\beta
k}^{\alpha}h_{ij}^{\alpha}h_{ij}^{\beta})+\sm{\sum}_{\beta} \bar{R}_{i\beta
j}^{\alpha}H^{\beta}h_{ij}^{\alpha}\nonumber\\
&&~~~+2\sm{\sum}_{kl}(\bar{R}_{ijk}^lh_{ij}^{\alpha}h_{kl}^{\alpha}
+\bar{R}_{kik}^lh_{lj}^{\alpha}h_{ij}^{\alpha})-\sm{\sum}_{k\beta}h_{ij}
^{\alpha}h_{jk} ^{\alpha}h_{ki} ^{\beta}H^{\beta}~\La{)}\nonumber\\
&&~~~+\sm{\sum}_{ij\alpha\beta}\La{[}\sm{\sum}_{k}(h_{ik} ^{\alpha}h_{jk}
^{\beta}-h_{ik} ^{\beta}h_{jk} ^{\alpha})\La{]}^2+\sm{\sum}_{
\alpha\beta}(\sm{\sum}_{ij}h_{ij}^{\alpha}h_{ij}^{\beta})^2 .
\end{eqnarray*}
We obtain a  Simons' type identity
\begin{eqnarray}
\Delta||B||^2&=&2||\nabla B||^2+\sm{\sum}_{ij\alpha}2
h_{ij}^{\alpha}H_{,ij}^{\alpha}-\sm{\sum}_{ij\alpha}2h_{ij}^{\alpha}
[\sm{\sum}_k(\bar{\nabla}_j\bar{R})_{kik}^{\alpha}+\sm{\sum}_k(\bar{\nabla}
_k\bar{R})_{ijk}^{\alpha}]\nonumber\\
&&+\sm{\sum}_{ij\alpha\beta}2\{\sm{\sum}_k(4\bar{R}_{\beta
ki}^{\alpha}h_{kj}^{\beta}h_{ij}^{\alpha}-\bar{R}_{k\beta
k}^{\alpha}h_{ij}^{\alpha}h_{ij}^{\beta})+\bar{R}_{i\beta
j}^{\alpha}H^{\beta}h_{ij}^{\alpha}\}\nonumber\\
&&+\sm{\sum}_{ijkl\alpha}4(\bar{R}_{ijk}^lh_{ij}^{\alpha}h_{kl}^{\alpha}
+\bar{R}_{kik}^lh_{lj}^{\alpha}h_{ij}^{\alpha})
-\sm{\sum}_{ijk\alpha\beta}2h_{ij}
^{\alpha}h_{jk} ^{\alpha}h_{ki} ^{\beta}H^{\beta}\nonumber\\
&&+2\sm{\sum}_{ij\alpha\beta}\La{(}\sm{\sum}_{k}(h_{ik} ^{\alpha}h_{jk}
^{\beta}-h_{ik} ^{\beta}h_{jk} ^{\alpha})\La{)}^2+2\sm{\sum}_{\alpha
\beta}(\sm{\sum}_{ij}h_{ij}^{\alpha}h_{ij}^{\beta})^2.\label{eq 2.4}
\end{eqnarray}
Notice that we use $||\cdot||$ to denote the absolute of the norm of
a time-like vector in $\oM$.

\section{$\Delta\cosh\theta$}
\renewcommand{\thesection}{\arabic{section}}
\renewcommand{\theequation}{\thesection.\arabic{equation}}
\setcounter{equation}{0}
\setcounter{theorem}{0}
\setcounter{corollary}{0}
\setcounter{proposition}{0}
\setcounter{lemma}{0}

In this section we shall compute the covariant derivative of a
pull-back of a parallel form in the ambient space by a spacelike
immersion $F:M^m\to\oM$. Let $\Omega $ be a
parallel $m$-form on $\oM$. For the orthonormal frame fields
$\{e_i, e_{\alpha}\}$ in section 2, $\Omega (e_1, \cdots, e_m)$ is
a function on $M$. As in \cite{w2,ls}, we shall compute the
Laplacian of $\Omega (e_1, \cdots, e_m)=\Omega _{1\cdots m}$ in
locally frame fields.
First we have
\begin{eqnarray}
(\nabla _kF^*\Omega)(e_{1}, \cdots, e_{m})&=&\sum_i\Omega(e_1,\ldots,
(\bar{\nabla} _ke_{i}-\nabla _ke_{i}), \ldots,e_{m} ) \nonumber\\
&=&\sum_i\Omega(e_1,\ldots, B(e_k,e_i),\ldots, e_m).\label{eq 3.1}
\end{eqnarray}
Differentiating (\ref{eq 3.1}) again gives
\begin{eqnarray*}
(\Delta F^*\Omega)(e_{1}, \cdots, e_{m}) &=&
\sum_i\Omega(e_1,\ldots,\sum_k(\nabla_{e_k}B)(e_k,e_i)
+(\bar{\nabla}_{e_k}B(e_k,e_i))^{\top},\ldots,e_m)\\
&&+\sum_k\sum_{j<i}\Omega( e_1, \ldots,
B(e_k,e_j),\ldots,B(e_k,e_i),\ldots,e_m)\\
&&+\sum_k\sum_{j>i}\Omega( e_1, \ldots,
B(e_k,e_i),\ldots,B(e_k,e_j),\ldots,e_m),\\
\end{eqnarray*}
where $\Delta F^*\Omega=\sum_k \nabla_k\nabla _kF^*\Omega-
\nabla_{\nabla_{e_k}e_k}F^*\Omega$ is the rough Laplacian.
Using the Codazzi's equation
$\sum_k\nabla^{\bot}_{e_k}B(e_k,e_i)=
\nabla^{\bot}_{e_i}H+(\bar{R}(e_k,e_i)e_k)^{\bot}$
and that
$$\sum_{ik}g((\bar{\nabla}{e_k}B(e_k,e_i))^{\top},e_i)=
\sum_{ik}-\bar{g}(B(e_k,e_i),B(e_k,e_i))= \|B\|^2,$$ we get in
components
\begin{equation}
(\Delta F^*\Omega )_{1\cdots m}= \Omega _{1\cdots
m}||B||^2+2\sum
_{\alpha <\beta, i<j}\Omega _{\alpha \beta ij}\hat{R}_{\beta ij}^{\alpha }
+\sum _{\alpha, i}\Omega _{\alpha i}H_{,i}^{\alpha }-\sum
_{\alpha, i, k}\Omega _{\alpha i}\bar{R}_{kik}^{\alpha }, \label{eq 3.2}
\end{equation}
where $ H_{,i}^{\alpha }=(\nabla^{\bot}_{e_i}H)^{\alpha}$, and
$$\hat{R}_{\beta ij}^{\alpha}=h_{ik}^{\alpha
}h_{jk}^{\beta}-h_{ik}^{\beta}h_{jk}^{\alpha},~~~~~~~~ \Omega
_{\alpha \beta ij}=\Omega (e_1, \cdots , e_{\alpha },\cdots,
e_{\beta},\cdots,e_m)$$
 with $e_{\alpha}$, $e_{\beta}$ occupying
the $i$-th and the $j$-th positions. The same meaning is for
$\Omega _{\alpha i}$.
\\[-1mm]

In the following we assume $\oM= \Sigma_1\times
\Sigma _2$ is a product of two Riemannian manifolds $(\Sigma_i, g_i)$
 of dimension $m$ and $n$,   with pseudo-Riemannian metric
$\bar{g}=g_1-g_2$. If we
denote by $\pi_i$ the projection from $T\oM$ onto $T\Sigma
_i$, then for any $X, Y\in T\oM$,
\begin{equation}\label{eq 3.3}
\bar{g}(X, Y)=g_1(\pi_1(X),\pi_1(Y))-g_2(\pi_2(X),\pi_2(Y)).
\end{equation}

Suppose $M$ is a spacelike graph of a smooth map $f: \Sigma
_1\rightarrow \Sigma _2$. For each $p\in \Sigma _1$
let $\lambda_1^2\geq
\lambda_2^2\geq \ldots\geq \lambda_m^2\geq 0$ be the eigenvalues
of $f^*g_2$. The spacelike
condition on $M$ means $\lambda_i^2<1$. By the classic Weyl's
perturbation theorem \cite{we}, ordering the eigenvalues in this
way, each $\lambda_i^2:\Sigma_1\rightarrow [0,1)$ is a continuous
locally Lipschitz function.
 For each $p$
let $ s=s(p)\in\{1\,\ldots, m\} $ be the rank of $f$ at $p$, that
is, $\lambda_s^2>0$ and $\lambda_{s+1}=\ldots=\lambda_m=0$. Then
$s\leq \min\{m,n\}$.

We take a $g_1$-orthonormal basis $\{a_i\}_{i=1, \cdots, m}$ of
$T_p\Sigma _1$ of eigenvectors of $f^*g_2$ with corresponding
eigenvalues $\lambda_i^2$. Set $a_{i+m}=df(a_i)/\|df(a_i)\|$ for
$i\leq s$. This constitutes an orthonormal system in
$T_{f(p)}\Sigma_2$, that we complete to give an orthonormal basis
 $\{a_{\alpha}\}_{\alpha =m+1, \cdots, m+n}$ for
$T_{f(p)}\Sigma _2$. Moreover, changing signs if necessary, we can
write
 $df(a_i)=-\lambda _{i\alpha}a_{\alpha}$,
where $\lambda _{i\alpha }=\delta_{\alpha,m+i}\lambda_i$
 meaning $=0$ if $i> s$,
or $\alpha> m+s$. Therefore
\begin{eqnarray}
e_i &=& \sm{\frac 1{\sqrt{1-\sum
_{\beta}\lambda _{i\beta}^2}}}(a_i+\sum _{\beta}\lambda
_{i\beta}a_{\beta})~~~~i=1, \cdots , m \label{eq 3.4}\\
e_{\alpha} &=&\sm{\frac 1{\sqrt{1-\sum _j\lambda _{j\alpha
}^2}}}(a_{\alpha}+\sum _j\lambda _{j\alpha}a_{j})~~~~\alpha
=m+1, \cdots , m+n \label{eq 3.5}
\end{eqnarray}
form an orthonormal basis for $T_pM$ and for $N_pM$ respectively,
with $e_i$ a direct one.

From now on we take $\Omega $ to be the volume form of $\Sigma
_1$, which is a parallel $m$-form on $\oM$. If $M$ is a
embedded $m$-submanifold such that for any $p\in M$,
and a basis $E_i$ of $T_pM$, the
quantity $\Omega (\pi_1(E_1), \cdots, \pi_1(E_m))$  is non-null
then $M$  is locally a graph,
for the
later implies $\pi_1\circ F:M\rightarrow \Sigma_1$ is a local
diffeomorphism. This means
$F(p)=(\phi(p),f(\phi(p)))$
where $\phi:M\rightarrow \Sigma_1$ is a local diffeomorphism,
and $F$ can be locally identified with the graph $\tilde{F}(p)=(p,f(p))$ up
to parameterization. The mean curvature of $F$ does not depend
on the parameterization, only on its image.
 We shall call graphs to all such parameterizations. 
Note that by Lemma 3.1 of \cite{alal}, if
$F$ is a spacelike submanifold with $M$ compact,
 $\phi:M\to \Sigma_1$ is a covering
map, and so it is surjective. Hence, 
this map $f:\Sigma_1\to \Sigma_2$,  when locally defined 
and $\Sigma_1$ compact, it is unique and 
globally defined.

Assume $M=\Gamma_f$. If $M$ is a spacelike graph, then taking the
orthonormal frame $e_i$ as in (\ref{eq 3.4})
$$\Omega _{1\cdots m}=\Omega (\pi_1(e_1), \cdots, \pi_1(e_m))=*F^*\Omega
={\frac 1{\sqrt{\prod _{i=1}^{m}(1-\lambda
_i^2)}}}={\frac{1}{\sqrt{det(g_1-f^*g_2)}}},
$$
where $*$ is the star operator in $M$. In this case this
 quantity is  $\geq 1$ (assuming the correct orientation)
 and is $\cosh\theta$. We can also describe $\cosh\theta$ as the
ratio between the volume elements of $(\Sigma_1,g_1)$ and of
$(\Sigma_1,g=g_1-f^*g_2)$.
If $M$ is compact, any other submanifold in a sufficiently small
neighbourhood of $M$ is also a spacelike graph. 
Now we compute
\begin{eqnarray}2\sum _{\alpha <\beta, i<j}\Omega _{\alpha\beta
ij}\hat{R}_{\beta ij}^{\alpha }=2\sum _{\alpha ,\beta , k, i<j
}\lambda _{i\alpha }\lambda
_{j\beta}(h_{ik}^{\alpha}h_{jk}^{\beta}-h_{ik}^{\beta}
h_{jk}^{\alpha})\cosh\theta. \label{eq 3.6}
\end{eqnarray}
As for the terms containing the curvatures of the ambient space,
we denote by $R_1$ and $R_2$ the curvature tensor of $\Sigma _1$
and $\Sigma _2$, respectively. We shall compute the curvatures
$\bar R$ of $\bar M$ in terms of $R_1$ and $R_2$. Now for the
tangent frame field $\{e_i\}$ (\ref{eq 3.4}) and normal frame field
$\{e_{\alpha}\}$ (\ref{eq 3.5}), since
$\bar{R}_{kik}^{\alpha}=\bar{R}_{\beta kik}\bar{g}^{\alpha
\beta}=-\bar{R}_{\alpha kik}$, we obtain
 \begin{eqnarray*}
&&-\bar{R}_{kik}^{\alpha}=\bar{R}(e_{\alpha},e_k, e_i, e_k)\nonumber \\
&=& R_1(\pi_1(e_{\alpha}),
\pi_1(e_k),\pi_1(e_i),\pi_1(e_k))-R_2(\pi_2(e_{\alpha}),
\pi_2(e_k),\pi_2(e_i),\pi_2(e_k))\\
&=&\! \frac {\sm{\sum}
_{l}\lambda_{l\alpha}R_1(a_l, a_k, a_i, a_k)\!-\!\!\!\sm{\sum}
_{\beta,\gamma,\delta }\!\lambda_{k\beta}\lambda_{i\gamma}\lambda
_{k\delta} R_2(a_{\alpha},a_{\beta}, a_{\gamma},
a_{\delta})}{\sm{\sqrt{(1-\sum_j\lambda _{j\alpha }^2)(1-\lambda
_{i}^2)}(1-\lambda _{k}^2)}}
.
\end{eqnarray*}
Consider for $i\neq j$ the two-planes $P_{ij}=span\{a_i,a_j\}$,
$P'_{ij}=span\{a_{m+i},a_{m+j}\}$.
Since $\lambda_{i\alpha}$ is diagonal, we have
\begin{eqnarray}
 &&-\sum_{\alpha, i, k}\Omega _{\alpha
i}\bar{R}_{kik}^{\alpha}\nonumber\\
~~~~~~~&=&\sum_{i,j}\frac{\cosh\theta\lambda_i^2}{(1-\lambda_i^2)
(1-\lambda_j^2)}\left(R_1(a_i, a_j, a_i, a_j)-
\lambda_j^2R_2(a_{m+i},a_{m+j}, a_{m+i},a_{m+j})\right)
\nonumber\\
&=& \cosh\theta\sum_{i, j\neq i}\left(\frac{ \lambda_i^2}{(1-\lambda
_i^2)}K_1(P_{ij})+\frac{ \lambda _i^2\lambda_j^2}
{(1-\lambda_i^2)(1-\lambda_j^2)} [K_1(P_{ij})-K_2(P'_{ij})]\right).~~~~~~~
\label{eq 3.7}
\end{eqnarray}
Inserting (\ref{eq 3.6}) and (\ref{eq 3.7}) into (\ref{eq 3.2})
we at last arrive at
\begin{eqnarray}
\lefteqn{\Delta \cosh\theta = \cosh\theta \La{\{}||B||^2+2\sum
_{k,i<j}\lambda _i\lambda _jh_{ik}^{m+i}h_{jk}^{m+j}-2\sum
_{k,i<j}\lambda
_i\lambda _jh_{ik}^{m+j}h_{jk}^{m+i}}\nonumber\\
~~~&&+\sum_{i}\La{(}\frac{\lambda _i^2}{(1-\lambda _i^2)}Ricci_1(a_i,a_i)
+\sum_{j\neq i}
{\frac{\lambda _i^2\lambda_j^2}{(1-\lambda_i^2)(1-\lambda_j^2)}}
\left[K_1(P_{ij})-K_2(P'_{ij})\right]\La{)}\La{\}}\nonumber\\
~~~&&~~~+\sum _{\alpha, i}\Omega _{\alpha i}H_{,i}^{\alpha}, \label{eq 3.8}
\end{eqnarray}
where we have used the fact that the Hodge star operator is
parallel.
Now by (\ref{eq 3.1}) we have
$
d\cosh\theta(e_k)
=\cosh\theta \sum _i\lambda _ih_{ik}^{m+i}$,
which implies
\begin{equation}\frac {|\nabla \cosh\theta|^2}{\cosh^2\theta
}=\sum _k(\sum _i\lambda _ih_{ik}^{m+i})^2=\sum
_{i,k}(\lambda_ih_{ik}^{m+i})^2+2\sum _{i<j,k}\lambda _i\lambda
_jh_{ik}^{m+i}h_{jk}^{m+j}. \label{eq 3.9}
\end{equation}
We shall calculate
\begin{equation}\Delta \ln(\cosh\theta)=\frac
{\cosh\theta \Delta(\cosh\theta)-|\nabla
\cosh\theta|^2}{\cosh^2\theta }. \label{eq 3.10}
\end{equation}
Plugging (\ref{eq 3.8}) and (\ref{eq 3.9}) into (\ref{eq 3.10}) we have
\begin{eqnarray}
\lefteqn{\Delta \ln(\cosh\theta)=||B||^2-\sum
_{i,k}\lambda_i^2(h_{ik}^{m+i})^2-2\sum _{k,i<j}\lambda
_i\lambda _jh_{ik}^{m+j}h_{jk}^{m+i}}\nonumber\\
~~~&&+\sum_{i}\La{(}\frac{\lambda _i^2}{(1-\lambda _i^2)}Ricci_1(a_i,a_i)
+\sum_{j\neq i}\frac{\lambda _i^2\lambda_j^2}{(1-\lambda_i^2)(1-\lambda_j^2)}
\left[K_1(P_{ij})-K_2(P'_{ij})\right]\La{)}\nonumber\\
~~~&&~~~+(\cosh\theta)^{-1}\sum _{\alpha , i}\Omega _{\alpha
i}H_{,i}^{\alpha }.\label{eq 3.11}
\end{eqnarray}

\section{Evolution equations}
\renewcommand{\thesection}{\arabic{section}}
\renewcommand{\theequation}{\thesection.\arabic{equation}}
\setcounter{equation}{0}
\setcounter{theorem}{0}
\setcounter{corollary}{0}
\setcounter{proposition}{0}
\setcounter{lemma}{0}

In this section, we shall compute the evolution equations of several
geometric quantities along the mean curvature flow (1.1).
We fix a point $(x_0,t_0)\in M\times [0, T)$ and
consider $(x,t)$ in a neighbourhood of $(x_0,t_0)$. We locally identify
$M_t=(M,g_t=F_t^*\bar{g}))$ with $(F_t(M),\bar{g}_{F_t(M)})$. We
take $e_{\alpha}(x,t)$ a local o.n. frame of $NM_t$ defined for
$(x,t)$ near $(x_0,t_0)$. Computations are easier considering a
fixed local coordinate chart on $M$. For any local coordinate
$\{x^i\}$ on $M$, we use the same notation as in section 3, but
 $g_{ij}(x,t)=g_t(\partial_i,\partial_j)=
\delta_{ij}$ may not hold everywhere, and $h_{ij}^{\alpha}$ the
components of the second fundamental form  $B(x,t)$ are with respect
to  $\partial_i$ and some  orthonormal frame $e_{\alpha}(x,t)$. That is
$$
{g}_{ij}(x,t)=\bar{g}(\frac {\partial F}{\partial x^i},\frac
{\partial F}{\partial x^j}),~~~~~ {h}^{\alpha}_{ij}(x,t) =-\bar{g}(
B(\partial_i,\partial_j), e_{\alpha}).
$$
Most of the following computations  are quite well known in the
literature ( see for instance \cite{cl,h,x2}), but for the sake
of simplicity, and since we are in the non-flat pseudo-Riemannian
setting and in higher codimension,  we reproduce them here
adapted to our case.
We define the tensor on $M$ (depending on $t$)
$$
 \mathcal{H}(X,Y) =-\bar{g}(B(X,Y),H)
$$
and $\mathcal{H}_{ij}=\mathcal{H}(\partial_i,\partial_j)$.
Using
the symmetry of the Hessian of
$F:M\times [0,T)\rightarrow \oM$ ( $M$ with
the initial metric $g_0$), and that $\nabla_{\frac{d}{dt}}\partial_i
=\nabla_{\partial_i}\frac{d}{dt}=0$,
we have
$$\bar{\nabla}_{H}\frac {\partial F}{\partial x^i}=
{\nabla}_{\frac{d}{dt}}{d F}(\partial_i)=
\bar{\nabla}_{\partial_i}\frac{d}{dt}F=\bar{\nabla}_{\partial_i}H.$$
Then, $\frac
d{dt}{g}_{ij}=\bar{g}(\bar{\nabla}_{\partial_i}H,
dF(\partial_j))+\bar{g}(dF(\partial_i),\bar{\nabla}_{\partial_j}H).$
It follows the induced metric evolves according to
\begin{equation}\begin{array}{l}
\frac d{dt}{g}_{ij}=\sm{\sum}_{\alpha}
2{H}^{\alpha}{h}_{ij}^{\alpha}=2\mathcal{H}_{ij}\\[2mm]
\frac d{dt}{g}^{ij}=-\sum_{kr\alpha}2{g}^{ik}{g}^{rj}
{H}^{\alpha}{h}_{kr}^{\alpha}=-\sum_{kr}2g^{ik}g^{rj}\mathcal{H}_{kr}
\end{array} \label{eq 4.1}
\end{equation}
for $(x,t)$ near $(x_0,t_0)$.
The volume element of $M_t$ is given by
  $d\mu_t= Vol_{M_t}= \sqrt{det[{g}_{ij}]}dx^{1\ldots m}$.
To compute the evolution equation for $d\mu_t$ and for the
second fundamental form  we will assume
the coordinate chart $x^i$ is normal at $x_0$ for the
metric $g_{t_0}$ with $\partial_i(x_0)=e_i(x_0)$ orthonormal frame.
Then at $(x_0,t_0)$,
${g}_{ij}=\delta_{ij}$.
The next
computations are at the point $(x_0,t_0)$. Using (\ref{eq 4.1})
$$ \frac{d}{dt}_{|t=t_0} d\mu_t=\frac{1}{2}\sm{\sum}_k
\frac{d {g}_{kk}}{dt}d\mu_{t_0}=
\|H\|^2 d\mu_{t_0}.$$
We also have
\begin{eqnarray*} \frac {d}{dt}{h}_{ij}^{\alpha}&=& -\frac
d{dt}\bar{g}(\bar{\nabla} _{\partial
_j}\partial _i, e_{\alpha})\\
&=&-\bar{g}(\bar{\nabla} _{{H}}\bar{\nabla} _{\partial _j}\partial _i,
e_{\alpha})-\bar{g}(\bar{\nabla} _{\partial _j}\partial _i, \bar{\nabla}
_{{H}}e_{\alpha})\\
&=&-\bar{g}(\bar{\nabla} _{\partial _j}\bar{\nabla} _{{H}}\partial
_i+\bar{R}({H},
\partial _j)\partial _i, e_{\alpha})
-\bar{g}(\bar{\nabla} _{\partial _j}\partial_i,\bar{\nabla} _{{H}}e_{\alpha})\\
&=&-\bar{g}(\bar{\nabla} _{\partial _j}\bar{\nabla} _{\partial _i}{H},
e_{\alpha})-\bar{g}(\bar{R}({H},
\partial _j)\partial _i, e_{\alpha})
-\bar{g}(\bar{\nabla} _{\partial _j}\partial
_i, \bar{\nabla} _{{H}}e_{\alpha}).
\end{eqnarray*}
 Set $H^{\alpha}_{,ij}=(\nabla^2_{e_i,e_j}H)^{\alpha}$.
Since at $(x_0,t_0)$, $\nabla_{\partial_i}\partial_j=0$ then
$$\begin{array}{l}
\bar{g}(\bar{\nabla} _{\partial _j}\bar{\nabla} _{\partial
_i}{H},e_{\alpha})=-H_{,ij}^{\alpha}- \sm{\sum}_{k \beta}
H^{\beta}h_{ik}^{\beta}h_{jk}^{\alpha}\\[1mm]
\bar{g}(\bar{\nabla} _{\partial _j}\partial
_i, \bar{\nabla}
_{{H}}e_{\alpha})=\sm{\sum}_{\beta}h_{ij}^{\beta}\bar{g}
(e_{\beta},\bar{\nabla}_{{H}}e_{\alpha}).\end{array}$$
Thus,
$$
\frac {d}{dt}{h}_{ij}^{\alpha}=H_{,ij}^{\alpha}
+\sm{\sum}_{k\beta}H^{\beta}h_{ik}^{\beta}h_{jk}^{\alpha}-\sm{\sum}_{\beta}
H^{\beta}\bar{R}_{i\alpha
j\beta}-h_{ij}^{\beta}\bar{g}(e_{\beta},\bar{\nabla}
_{{H}}e_{\alpha})
$$
and using $\sum_{\alpha\beta}h^{\alpha}_{ij}h^{\beta}_{ij}\bar{g}(e_{\beta},
\bar{\nabla}_{H}e_{\alpha})=0$,
 we have at $(x_0,t_0)$
\begin{eqnarray*}\frac d{dt}||B||^2&=& \sm{\sum}_{ijls}\frac d{dt}(
{g}^{il}{g}^{js}{h}_{ij}{h}_{ls})\\
&=&\sm{\sum}_{ijk\alpha\beta}-4H^{\beta
}h_{ij}^{\beta}h_{kj}^{\alpha}h_{ki}^{\alpha}+
\sm{\sum}_{ij\alpha}2h_{ij}^{\alpha}\frac
d{dt}h_{ij}^{\alpha}\\
&=&\sm{\sum}_{ij\alpha}2h_{ij}^{\alpha}H_{,ij}^{\alpha}
- \sm{\sum}_{ij\alpha\beta}(2H^{\beta
}h_{ij}^{\beta}h_{kj}^{\alpha}h_{ki}^{\alpha}+
2H^{\beta}\bar{R}_{i\alpha j\beta}h_{ij}^{\alpha}).
\end{eqnarray*}
\noindent
Combining the above equation with the Simon's type identity (\ref{eq 2.4}) we
arrive last to the evolution equation of the squared norm of the second
fundamental form as stated in next proposition. A similar
 computation can  be
done to $\|H\|^2$. Therefore
\begin{proposition}  \label{proposition 4.1}
Let $F: M\times [0, T)\rightarrow \oM$
be an $m$-dimensional mean curvature flow of a spacelike
submanifold in a pseudo-Riemannian manifold $\oM$. Then the
following evolution equations hold at $(x_0,t_0)$
\begin{eqnarray*}
\frac{d}{dt}d\mu_t &=&\|H\|^2d\mu_t\\
\frac{d}{dt}\|H\|^2
&=& \Delta\|H\|^2-2\|\nabla^{\bot}H\|^2-4\|\mathcal{H}\|^2
-2trace_g\bar{R}(dF(\cdot),H,dF(\cdot),H) \\[2mm]
\frac d{dt}||B||^2
&=&\Delta ||B||^2-2||\nabla B||^2+{\sum}_{ij\alpha}2h_{ij}^{\alpha}\La{(}
\sm{\sum}_k(\bar{\nabla}
_j\bar{R})_{kik}^{\alpha}+(\bar{\nabla}
_k\bar{R})_{ijk}^{\alpha}\La{)}\nonumber\\
&-&\!2\left({\sm{\sum}}_{ijk\alpha\beta}(4\bar{R}_{\beta
ki}^{\alpha}h_{kj}^{\beta}h_{ij}^{\alpha}-\bar{R}_{k\beta
k}^{\alpha}h_{ij}^{\alpha}h_{ij}^{\beta})
+{\sm{\sum}}_{ijkl\alpha}2(\bar{R}_{ijk}^lh_{ij}^{\alpha}h_{kl}^{\alpha}
+\bar{R}_{kik}^lh_{lj}^{\alpha}h_{ij}^{\alpha})\!\!\right)\nonumber\\
&&-2\sm{\sum}_{ij\alpha\beta}\La{(}\sm{\sum}_{k}(h_{ik} ^{\alpha}h_{jk}
^{\beta}-h_{ik} ^{\beta}h_{jk} ^{\alpha})\La{)}^2-2\sm{\sum}_{\alpha,
\beta}(\sm{\sum}_{ij}h_{ij}^{\alpha}h_{ij}^{\beta})^2.\\
\end{eqnarray*}
\end{proposition}

Next we compute the evolution of the pull-back of a parallel
$m$-form on $\oM$. Let $\Omega $ be a parallel $m$-form on
$\oM$. Then the restriction of $F^*\Omega $  satisfies the
following evolution equation at $(x_0,t_0)$
\begin{eqnarray*}
\lefteqn{\frac d{dt} (F^*\Omega (\partial _1, \cdots, \partial _m))=\frac
d{dt}(\Omega (F_*\partial _1, \cdots, F_*\partial _m))}\\
&=&\sum_i\Omega ( \partial _1, \cdots,
\nabla ^{\bot}_{\partial_i}{H}, \cdots,
\partial _m)+ \Omega (\partial _1, \cdots,
-A_{ {H}}{\partial_i},  \cdots,\partial _m)\\
&=&\sum_{\alpha, i}\Omega ( \partial _1, \cdots,
e_{\alpha}, \cdots,\partial _m)H_{,i}^{\alpha}
+\Omega (\partial _1, \cdots,
\partial _m)\sum _{\alpha}(H^{\alpha})^2\\
&=&\sum_{\alpha, i} \Omega _{\alpha i}H_{,i}^{\alpha} +\cosh\theta
||{H}||^2,
\end{eqnarray*}
 On the other hand we have
\begin{eqnarray}
\frac d{dt}\Omega _{1\cdots m}&=& \frac d{dt}(\Omega (\partial_ 1,
\cdots, \partial _m)\frac 1{\sqrt{g}})\nonumber\\
&=&\frac 1{\sqrt{g}}\frac d{dt}\Omega (\partial_ 1, \cdots,
\partial _m)-\frac 1{\sqrt{g}}||{H}||^2\Omega (\partial_ 1,
\cdots, \partial _m)\nonumber\\
&=& \sum_{\alpha, i} \Omega_{\alpha i}H_{,i}^{\alpha}. \label{eq
4.2}
\end{eqnarray}
Combining with equation (3.2) we get the parabolic equation
satisfied by $\Omega _{1\cdots m}$:
\begin{proposition} \label{proposition 4.2}
Let $M_t$ be an $m$-dimensional spacelike mean curvature
flow in a pseudo-Riemannian manifold $\oM$ and $\Omega $  a
parallel $m$-form on $\oM$. Then we have the following
evolution equation at $(x_0,t_0)$
\begin{eqnarray*}
\frac d{dt}\Omega _{1\cdots m}&=&\Delta \Omega_{1\cdots m}-
\Omega _{1\cdots m}||B||^2\nonumber\\
&&-2\sum _{\alpha<\beta, i<j}\Omega _{\alpha \beta
ij}\hat{R}_{\beta ij}^{\alpha}+\sum _{\alpha, i, k}\Omega _{\alpha
i}\bar{R}_{kik}^{\alpha}.
\end{eqnarray*}
\end{proposition}
If $M$ is a graph of $f:\Sigma _1\rightarrow \Sigma_2$, and
$\Omega $ is the volume form of $\Sigma _1$, then, since $\Sigma_1$
is compact, for sufficiently small $t$, $M_t$ is a spacelike graph,
and so $\cos\theta_t$  is defined and
we have the
evolution equation for $\cosh\theta$ by inserting (\ref{eq 3.11}) into
(\ref{eq 4.2})
\begin{proposition}  \label{proposition 4.3}
Let $F_0:M\rightarrow \oM$ be an immersion such that $M_0$ is  a
spacelike graph over $\Sigma _1$. If each $M_t$ is  a graph
$\Gamma_{f_t}$ of a map $f_t:\Sigma _1\rightarrow \Sigma_2$
along the mean curvature flow of $F_0$ for $t\in
[0, T')$, $T'\leq T$,
then $\cosh\theta$ satisfies the following equation, using the frames
(\ref{eq 3.4})
and (\ref{eq 3.5})
\begin{eqnarray}
\lefteqn{\frac d{dt}\ln(\cosh\theta) =\Delta \ln(\cosh\theta)~+}\nonumber\\
&&~~~~~-
\La{\{}||B||^2-\sum _{k,i}\lambda _i^2(h_{ik}^{m+i})^2-2\sum
_{k,i<j}\lambda
_i\lambda _jh_{ik}^{m+j}h_{jk}^{m+i}\La{\}} \label{eq 4.3}\\
\!\!&&-\sum_{i}\La{(}\frac{\lambda _i^2}{(1-\lambda
_i^2)}Ricci_1(a_i,a_i) +\sum_{j\neq i}\frac{\lambda
_i^2\lambda_j^2}{(1-\lambda_i^2)(1-\lambda_j^2)}
\left[K_1
(P_{ij})-K_2(P'_{ij})\right]\La{)}~~~~~~~~ \label{eq 4.4}
\end{eqnarray}
\end{proposition}

\section{Short-time existence}
\renewcommand{\thesection}{\arabic{section}}
\renewcommand{\theequation}{\thesection.\arabic{equation}}
\setcounter{equation}{0}
\setcounter{theorem}{0}
\setcounter{corollary}{0}
\setcounter{proposition}{0}
\setcounter{lemma}{0}

In this section we give the  proof of the first part
of Theorem 1.1(1).

Let $x^i$ be a coordinate chart of $\Sigma_1$ on a neighbourhood of
$p_0\in \Sigma_1$
and $y^{\alpha}$ a coordinate chart of
$\Sigma_2$ on a neighbourhood of $q_0= f_0(p_0)$.
Note that $x^i$ is identified with a coordinate chart in $M_0$
 as in section 4. In coordinates (1.1) means
\begin{equation}\label{eq 5.1a}
 \sum_{ij}{g}^{ij}\left(\frac{
\partial F_t^a}{\partial x_i\partial x_j}-
\sum_k\Gamma^k_{ij}\frac{\partial F^a_t}{\partial x_k}\right)=
-\bar{G}(x,t)^a +\frac{d}{dt}F^a_t,~~~a=1,\ldots,m+n
\end{equation}
where
$\bar{G}(x,t)^a = \sum_{ijbc}{g}^{ij}
(\bar{\Gamma}^{a}_{bc}\circ F_t)\, \frac{\partial F^{b}_t}{\partial
x^i} \frac{\partial F^{c}_t}{ \partial x^j}$,
and $\Gamma^k_{ij}$ are the Christoffel symbols for the induced
Riemannian metric $g_t$ of $M$ (that depends on the second
derivatives of $F_t$, what makes the system to be not strictly
parabolic)
 and $\bar{\Gamma}^{a}_{bc}$ the ones of $\bar{M}$,
in the coordinates charts $x^i$ and $w^{a}=(x^i,y^{\alpha})$
respectively.

Since $F=F_0$ is  a spacelike graph
 of $f=f_0:\Sigma_1\rightarrow \Sigma_2$
we recall that in \cite{sal}, following \cite{sal0},
 we have proved that for $X,Y\in T_p\Sigma_1$,
\begin{eqnarray}
B(X,Y) &=& (\nabla^1_XY-\nabla_XY,df(\nabla^1_XY-\nabla_XY))
+(0,Hess \, f(X,Y)) \nonumber \\
&=&(0,Hess\, f(X,Y))^{\bot}, \label{eq 5.1}\\[1mm]
H &=& (Z,df(Z))+(0,W)=(0,W)^{\bot}\nonumber
\end{eqnarray}
where
$Hess\, f$ is the Hessian of $f$ with respect to the
Levi-Civita  connections
$\nabla^i$ of $(\Sigma_i, g_i)$,
$W=trace_{g}Hess \, f$, and $Z$ is the vector field on $\Sigma_1$ defined by
$g(Z,X)=g_2(W,df(X))$.
From the above expression  of $B$ we have observed in
\cite{ls} that
$\Gamma_f$ is a totally geodesic submanifold of
$\oM$ iff $f:(\Sigma_1,g_1)\to (\Sigma_2,g_2)$ is a totally
geodesic map. To see this, we have from (\ref{eq 5.1}) and using the
frames (\ref{eq 3.4}) and (\ref{eq 3.5}), for $u\in T_q(\Sigma_2)$,
$u=\sum_{\alpha}u^{\alpha}a_{\alpha}$,
$$(0,u)^{\bot}=\La{(}\sum_i \frac{\lambda_i}{1-\lambda_i^2}u^{i+m}a_i,
\sum_i \frac{1}{1-\lambda_i^2}u^{i+m}a_{i+m} +\sum_{\alpha>2m}
u^{\alpha}a_{\alpha}\La{)}.$$
Using these frames we also see that
$\pi_1:TM_p\to T_p\Sigma_1$ and $\pi_2:NM_p\to T_{f(p)}\Sigma_2$
define  isomorphisms. From (\ref{eq 5.1}),
\begin{equation}
\begin{array}{l}
\sum_i \frac{1}{(1-\lambda_i^2)}
\|Hess\, f(\partial_k, \partial_j)^{i+m}\|^2_2 +
\sum_{\alpha>m}\|Hess\, f(\partial_k, \partial_j)^{\alpha}\|^2_2
=\|B(\partial_k, \partial_j)\|^2 \\[3mm]
\sum_k((\Gamma^1)^k_{ij}-\Gamma^k_{ij})\partial_k
=\pi_1(B(\partial_i, \partial_j))=\sum_k \frac{\lambda_k}
{(1-\lambda_k^2)}Hess\, f(\partial_i, \partial_j)^{k+m}a_{k}.
\end{array} \label{eq 5.2}
\end{equation}
We also note that if $\lambda_i^2<1-\delta$ for all $i$,
the Riemannian  metric $\hat{g}$ on $T_{(p,f(p))}\oM$
defined by declaring an orthonormal basis
 $e_i, e_{\alpha}$ given by (\ref{eq 3.4}) (\ref{eq 3.5}), that is,
 $\hat{g}= \bar{g}_{|{T\Gamma_{f}}}-\bar{g}_{|{T\Gamma_{f}}^{\bot}}$, 
 is equivalent to the Riemannian metric $\bar{g}_+=g_1+g_2$
of $\oM$  with
$ c(\delta)\bar{g}_+\leq \hat{g}\leq c'(\delta)\bar{g}_+$ along
$\Gamma_f$, 
where $c(\delta),c'(\delta)$ are positive constants that only
depend on $\delta$.
\\[4mm]
If  $F_t:M_0=\Sigma_1\rightarrow \bar{M}$
is   a  graph,
$ F_t(p)=(\phi_t(p),f_t(\phi_t(p)) )$,
where $\phi_t:\Sigma_1\rightarrow \Sigma_1$ is given by
$\phi_t(p)=\pi_1(F_t(p))$ and  satisfies $\phi_0=Id$,
then (1.1) means
$\frac{d\phi_t}{dt} = Z_t$, $\frac{df_t}{dt} =W_t$,
$\phi_0=Id,$ and $f_{t=0}=f_0$.
In particular, $f_t$ satisfies  the evolution equation
$$
\left\{\begin{array}{l}
\frac{df}{dt}= W_t= trace_{g_t}Hess\, f_t,\\[1mm]
f_{t=0}=f_0,
\end{array}\right.
$$
where the Hessian is w.r.t the initial metric $g_1$ of $\Sigma_1$
and  the trace with respect to the graph metric $g_t=g_1-f_t^*g_2$
of $\Sigma_1$.
 This  system  is strictly parabolic.

Now we assume $F_t$ satisfies (1.1).
We identify $M_0=F_0(M)$ with
the graph $\Gamma_{f_0}:\Sigma_1\rightarrow \bar{M}$. We also remark that
using the trick of DeTurck (see page 17 of \cite{zhu}), as in the case of
hypersurfaces in a Euclidean space, by reparameterizing
$F$ as $\hat{F}(p,t)=F(\rho_t(p),t)$ where $\rho_t:\Sigma_1\to \Sigma_1$
is a convenient (local) diffeomorphism,
(1.1) is equivalent to a system of strictly parabolic
equations.
For existence of short time solutions one can follow the approach in 
\cite{es,jost} of
isometrically embedding $\Sigma_i$ into Euclidean spaces
$\R^{N_i}$, but considering the Riemannian structures
on $\oM$ and $\R^{N_1+N_2}$,  and linearizing the above
parabolic system to prove existence of a local solution.
As we will see in next section,
the tensor fields involved, namely
  $\bar{\nabla}^kB$
are  bounded both for the pseudo-Riemannian and
the Riemannian structure of $\oM$. Since $\Sigma_1$ is compact one has:
\begin{proposition} \label{proposition 5.1}
A unique smooth solution of (1.1) with initial condition $F_0$ a
spacelike graphic submanifold exists in a maximal time interval $[0,
T)$ for some $T>0$.
\end{proposition}
\noindent
Let $T'\leq T$ such that for all $t<T'$,  $M_t$ is an entire spacelike
graph $\Gamma_{f_t}$, and
$\cosh\theta$ is bounded from above, that is
 $\cosh\theta =1/\sqrt{\prod_{i=1}^{m}(1-\lambda
_i^2)}<\Lambda$ for a constant $\Lambda>1$. This is equivalent to
$\lambda _i^2\le 1-\delta$ for some $\delta
>0$ and any $1\le i\le m$.
 We set
\begin{equation} \label{eq 5.3}
\eta_t:=\max_{M_t}\cosh\theta.
\end{equation}
Now we prove part (1) of theorem 1.1
\begin{proposition} \label{proposition 5.2} $T'=T$, that is
$\cosh\theta$ has a finite upper bound,  the
evolving submanifold $M_t$ remains a spacelike graph of a map $f_t:\Sigma
_1\rightarrow \Sigma_2$ whenever the flow (1.1) exists.
In particular $f_t^*g_2$ and $\|df_t\|^2$
(norm  with respect to the initial metric
$g_1$ of $\Sigma_1$) are uniformly bounded and
the Riemannian metrics $g_t$ on $\Sigma_1$ are
uniformly equivalent. Moreover
$$
 \int_0^T\sup_{\Sigma_1}\|H_t\|^2dt<c_0
$$
for some constant $c_0 >0$.
\end{proposition}
\noindent
{\em Proof.}
 Let $t<T'$.
Note that  $\lambda
_i\lambda _j<1-\delta $ for any $i$ and $j$, and  $\lambda
_i=0$ for $i>\min (m,n)$. For the second fundamental form, we have
$$||B||^2\ge \sum _{i,k, j}(h_{ik}^{m+j})^2=\sum
_{i<j,k}[(h_{ik}^{m+j})^2+(h_{jk}^{m+i})^2]+\sum _{i,
k}(h_{ik}^{m+i})^2,$$ where we keep in mind that $h_{ik}^{m+j}=0$
when $m+j>m+n$.
Therefore we can estimate the terms in the bracket of (\ref{eq 4.3})
\begin{eqnarray}
\quad&&||B||^2-\sum _{i,k}\lambda_i^2(h_{ik}^{m+i})^2-2\sum
_{k,i<j}\lambda _i\lambda
_jh_{ik}^{m+j}h_{jk}^{m+i}\nonumber\\
&\ge&\delta ||B||^2+(1-\delta)\left\{\sum
_{i<j,k}[(h_{ik}^{m+j})^2+(h_{jk}^{m+i})^2]+\sum _{i,
k}(h_{ik}^{m+i})^2\right\}\nonumber\\
&&-(1-\delta)\sum _{i,k}(h_{ik}^{m+i})^2-2(1-\delta)\sum
_{k,i<j}|h_{ik}^{m+j}||h_{jk}^{m+i}|\nonumber\\
&\ge& \delta ||B||^2. \label{eq 5.4}
\end{eqnarray}
On the other hand, since  $Ricci_1\geq 0$
and  $K_1(p)\ge K_2(q) $,
(\ref{eq 4.4}) is nonpositive.
Thus by Proposition 4.3, $\ln (\cosh\theta)$ satisfies the
differential inequality for all $t<T'$
$$
\frac d{dt}\ln(\cosh\theta) \le \Delta \ln(\cosh\theta)-\delta
||B||^2\leq \Delta \ln(\cosh\theta).
$$
According to the maximal principle for parabolic equations,
we have for $s>t$
\begin{equation}
\eta_s\leq \eta_t\leq \eta_0. \label{eq 5.5}
\end{equation}
Assume $T'<T$. Then $F_{T'}$ is defined, and from
$F_t^*\Omega_{1\ldots m}\geq 1$, for all $t<T'$ we obtain
the same for $t=T'$. Then $F_{T'}$ is a graph of a map
$f_{T'}$. From (\ref{eq 5.5}) we have for all $T'>t\geq 0$,
$\lambda^2_i(t)<1$ and
$$1-\lambda_i^2(t)\geq \prod_i(1-\lambda_i^2(t))\geq 
\frac{1}{\eta_t^2}\geq \frac{1}{\eta_0^2},$$ and so the same also
holds for $t=T'$, what proves that $f_{T'}$ defines a spacelike
graph $F_{T'}$. Thus, we may take  $T'=T$. Therefore,
$F_t(p)=(\phi_t(p),f_t(\phi_t(p)))$ with $\phi_t:\Sigma_1\to
\Sigma_1$ a covering map  homotopic to the identity, and so of degree
one, and necessarily orientation preserving local diffeomorphism, what implies
$\phi_t$ is also 1-1. 
Now it follows
that $g_t$ are
all uniformly equivalent and  $d\mu_t$ is uniformly
bounded in $\Sigma_1\times [0,T)$. From proposition \ref{proposition 4.1},
$d\mu_t=e^{\int_0^t\|H_s\|^2ds}d\mu_0$, and so $
\sup_{\Sigma_1}\int_0^{T}\|H_s\|^2ds<c_0$, for some constant
$c_0>0$. If we take $p_s\in \Sigma_1$ such that
$\|H_s\|(p_s)=\max_{\Sigma_1} \|H_s\|$, then we have
$\int_0^{T}\|H_s\|^2(p_s)ds\leq
\sup_{\Sigma_1}\int_0^{T}\|H_s\|^2ds<c_0$.
\qed\\[4mm]
We will need the following lemmas:
\begin{lemma} \label{lemma 5.1} \cite{eh}
Let $f$ be a function on $M\times [0, T_1]$ satisfying
$$(\frac d{dt}-\Delta)f\le -a^2f^2+b^2$$
for some constants $a, b \in \R$. Then we have
$f\le \frac ba+\frac 1{a^2t}$
everywhere on $M\times (0, T_1]$.
\end{lemma}
\begin{lemma} \label{lemma5.2} \cite{g2}
Let $\Sigma_1$ be a compact Riemannian manifold and
$f\in C^1(\Sigma_1 \times J)$ where $J$ is an open interval, then
$f_{\max}(t)=\max_{\Sigma_1}f(\cdot, t)$
 is  Lipschitz  continuous and there holds
a.e.
$\frac{d f_{\max}}{dt}(t)=\frac{d f}{d t}
(x_t,t),$
 where $x_t\in \Sigma_1$ is a point which the
maximum is attained.
\end{lemma}
 The Riemannian metrics $\hat{g}=\bar{g}_{|TM_t}-\bar{g}_{|TM_t^{\bot}}$
on $\oM$ defined along
the flow are uniformly equivalent to $\bar{g}_+$.
 Therefore, if $U$ is a vector
field of $\oM$ defined along the flow, and if $U$ is normal
or tangential to the flow, that is
$U=U^{\bot}$ or $U=U^{\top}$,
then $U$  is uniformly $\bar{g}$-bounded
if and only if it  is uniformly
$\bar{g}_+$-bounded. 
Hence, any vector field $U$  with $U^{\top}$ and
$U^{\bot}$  uniformly $\bar{g}$-bounded,   is also uniformly
$\bar{g}_+$-bounded.
\begin{proposition} \label{proposition 5.3}
$||B||$, $\|H\|$,   $\|\nabla^{k} B\|$,
and $\|\nabla^{k} H\|$, for all $k$,
 are uniformly bounded as long as the solution exists. Furthermore
 $\|B\|_{\bar{g}_+}$, $\|\bar{\nabla}^{k}B\|_{\bar{g}_+}$ and 
$\|\bar{\nabla}^kH\|_{\bar{g}_+}$
are uniformly bounded as well.
\end{proposition}
\noindent
{\em Proof.}
In our case  $\Sigma _i$ are of
bounded curvature tensors and their  covariant derivative.
By proposition \ref{proposition 5.2},   $e_i$ and $e_{\alpha}$ given by
(\ref{eq 3.4}) and (\ref{eq 3.5}) are uniformly
bounded. Thus, the terms in the expressions
in Proposition 4.1  involving the curvature tensor $\bar{R}$
are bounded.
 It follows from
Proposition 4.1, for some constants $c_1, c_2, c_3\geq 0$,
\begin{eqnarray} \label{eq 5.6}
\frac d{dt}||B||^2\le \Delta ||B||^2+c_1||B||+c_2||B||^2-\frac
2n||B||^4 \le \Delta ||B||^2-\frac 1n||B||^4+c_3,~~~
\end{eqnarray}
where we have used some elementary geometric-arithmetic  inequality
$$\sum _{\alpha, \beta}(\sum_{ij}h_{ij}^{\alpha}h_{ij}^{\beta})^2
\ge \sum _{\alpha}
(\sum _{ij}(h_{ij}^{\alpha})^2)^2\ge \frac 1n (\sum
_{i,j,\alpha}(h_{ij}^{\alpha})^2)^2=\frac 1n ||B||^4.$$
 Then applying Lemma 5.1 to (\ref{eq 5.6}), $\|B\|$ is uniformly bounded.
Then we proceed as  in \cite{e,eh,hamil,h,h3,h4,s} to prove boundedness
of $\|\nabla^kB\|$, using an interpolation formula for tensors.
We note that all terms including the curvatures of ambient space are of
lower orders of the second fundamental form than the main part.
Since $\nabla^kH =trace \nabla^kB$ 
 we obtain uniform boundedness
for $\nabla^kH$ for all $k\geq 0$.
$\nabla^{k}B$ is the $k$-order covariant derivative
of $B$ using the normal connection $\nabla^{\bot}$ of $NM$.
Now we prove boundedness of $\bar{\nabla}^kB$ in $T\oM$ for
the Riemannian structure $\bar{g}_+=g_1+g_2$.
For $k=0$, by the first equation of (\ref{eq 5.2}),
uniform boundedness of $\|B\|$
and of $\lambda_i$ implies
uniform boundedness of $\| Hess f\|_2$. It follows now
by second equation  of (\ref{eq 5.2}) we get uniform boundedness
of $\|\pi_1(B)\|_1$ and of $\Gamma^k_{ij}$.
In particular $\|\pi_2(B)\|_2$ is also bounded, and
this proves uniformly boundedness of $\|B\|_{\bar{g}_+}$. For $k=1$,
we have $\bar{\nabla}_{\partial_s}B(\partial_i,\partial_j)=
(\bar{\nabla}_{\partial_s}B(\partial_i,\partial_j))^{\top}+
{\nabla}_{\partial_s}B(\partial_i,\partial_j)$
and so
$$\hat{g} (\bar{\nabla}_{\partial_s}B(\partial_i,\partial_j),
 \bar{\nabla}_{\partial_s} B(\partial_i,\partial_j))=
\sum_{\alpha\beta lr}
{g}^{lr}h^{\alpha}_{ij}h^{\alpha}_{sl}h^{\beta}_{ij}h^{\beta}_{sr}
 +\|{\nabla}_{\partial_s}B(\partial_i,\partial_j)\|^2
$$
that is,
$\|\bar{\nabla}B\|^2_{\hat{g}}=\|(\bar{\nabla}B)^{\top}\|^2
+\|\nabla B\|^2 \leq c_{22}\|B\|^4+
\|\nabla B\|^2,$
with $c_{22}>0$ a constant not depending on $t$. Thus, $\bar{\nabla}B$
is uniformly $\hat{g}$-bounded and so $\bar{g}_+$ uniformly bounded.
Inductively we obtain the same for higher order derivatives.
\qed

\begin{corollary} \label{corollary 5.1}
$f^* g_2$ and  $\Gamma^{r}_{ij}$ and
their derivatives are uniformly bounded.
\end{corollary}
\noindent
{\em Proof.}
Uniform boundedness of $\|\bar{\nabla}^{r}B\|_{g_+}$ implies by
(\ref{eq 5.1})
and proposition \ref{proposition 5.2},
inductively on $r\geq 0$, uniform boundedness of
$\nabla^{r}_{\partial_s}\Gamma^k_{ij}$ and of
$\|\nabla^r Hess f\|$. Since
$$\nabla_{\partial_s}f^*g_2(\partial_i, \partial_j)
=g_2(Hess f(\partial_s,\partial_i),df( \partial_j))+
g_2(df( \partial_i),Hess f(\partial_s, \partial_j))$$
we obtain  the uniform boundedness of $f^*g_2$ and its  derivatives.
\qed
\section{Long-time existence and convergence}
\renewcommand{\thesection}{\arabic{section}}
\renewcommand{\theequation}{\thesection.\arabic{equation}}
\setcounter{equation}{0}
\setcounter{theorem}{0}
\setcounter{corollary}{0}
\setcounter{proposition}{0}
\setcounter{lemma}{0}

It is well known
that if $\|B\|$ is uniformly bounded then
the mean curvature flow exists for all time.
This is well known  for  hypersurfaces,
as in the above references,
and for the case of flat ambient space. For non flat space
and higher codimension in the Riemannian case see
 \cite{chli,cl,s}). This  holds as well in our setting.
For the sake of completeness
 we will apply Schauder theory for elliptic
 systems  to prove long-time existence and
a condition for the convergence of the flow  at infinity.

We are assuming $(\Sigma_1, g_1)$ compact and  $(\Sigma_2,g_2)$
 complete. In this section we  consider $\oM$ with the Riemannian metric
$\bar{g}_+=g_1+g_2$ and we may embed isometrically $\Sigma_i$ into an
Euclidean space $\R^{N_i}$,
and consider $\R^{N}$, $N=N_1+N_2$ with its
Euclidean structure. Note that $\bar{M}$ is a closed subset
of $\R^{N}$, and so if $K\subset \R^{N}$ is a
compact set for the Euclidean topology, then
$K\cap \oM$ is a compact set for $(\oM, \bar{g}_+)$.
Then we follow as in \cite{jost}.
For each $0\leq \sigma<1$ and $k<+\infty$ integer, the
spaces $C^{k+\sigma}(\Sigma_1,\oM)$  are endowed
 with the usual $C^{k+\sigma}$-H\"{o}lder norms 
(well defined up to equivalence using coordinates charts).
$C^{k+\sigma}(\Sigma_1,\oM)$ is a Banach manifold
 with tangent space
$C^{k+\sigma}(\Sigma_1,F^{-1}T\oM)$
at $F\in C^{k+\sigma}(\Sigma_1,\oM)$.
These spaces can be seen
as closed subsets of $C^{k+\sigma}(\Sigma_1,\R^{N})$.
We consider $F_t$ a solution of (1.1)
with initial condition $F_0=\Gamma_{f_0}$.
Then $F_t$ satisfies the parabolic system \ref{eq 5.1a}. Set
$a_{ij}={g}^{ij}(x)$, and
$b_k=\sum_{ij}{g}^{ij}\Gamma^k_{ij}(x)$.
From the computations in section 4, propositions 5.2 and
5.3,  and corollary \ref{corollary 5.1},
  we have that $a_{ij}$, $b_i$, $\frac{d}{d x_k}a_{ij}$,
 $\frac{d}{d x_k}b_{ij}$
(as well $\frac{d}{dt}a_{ij}$ and $\frac{d}{dt}b_{ij}$ as one
can easily check)
are uniformly bounded in $\Sigma_1$, and so
$a_{ij}$ and $b_k$
are uniformly $C^{1}$-bounded in
$\Sigma_1$.
Note that if a vector field $V$ of $\oM$ along the flow, is
${C^{1}(\Sigma_1,F^{-1}T\oM)}$ -uniformly bounded
for $\bar{g}_+$, then it is also in $\R^N$.
The same conclusion for the higher
order H\"{o}lder norms.
Note also from proposition 5.2,  $\|df_t\|^2$ is uniformly bounded,
and so $\|\frac{\partial F_t^a}{\partial x_i}\|_{\bar{g}_+}$ is
uniformly bounded. From
uniform boundedness of $\|B\|$ and of $\Gamma^{k}_{ij}$ established
in corollary \ref{corollary 5.1},
and that $\frac{\partial^2F^a}{\partial x_i,\partial x_j}=
B(\partial_i,\partial_j)^a +\sum_k\Gamma^{k}_{ij}(x)\frac{\partial F^a}
{\partial x_k}$,
we have uniform boundedness of
$\|\frac{\partial F_t^a}{\partial x_i}\|_{C^{\sigma}(\Sigma_1,\R^N)}$.
By  elliptic Schauder theory
 (see \cite{jost} p. 79 ),
we have that a solution $F_t$ of (5.1)
satisfies for each $t$
$$\begin{array}{l}\|F(\cdot, t)\|_{C^{1+\sigma}(\Sigma_1,\oM)}\leq
c_{4}\,( \|\bar{G}(\cdot,t)\|_{L^{\infty}(\Sigma_1,\R^N)}+
\|H\|_{L^{\infty}(\Sigma_1,\R^N)})\\
\|F(\cdot, t)\|_{C^{2+\sigma}(\Sigma_1,\oM)}\leq
\, c_5 (\|\bar G(\cdot,t)\|_{
C^{\sigma}(\Sigma_1,\R^N)}+ \|H\|_{C^{\sigma}(\Sigma_1,\R^N)}
),\end{array}$$
Here $c_i$, $i=0,1,\ldots $ are positive constants not depending on $t$.
Since $\|H\|_{C^{k+\sigma}(\Sigma_1,\R^N)}$ is uniformly bounded, we have
\begin{proposition} \label{proposition 6.1}
 Let  $F_t$ be a solution of (1.1) for
$t\in[0,T)$. If
 $\|F(\cdot, t)\|_{L^{\infty}(\Sigma_1,\oM)}$ is
uniformly bounded then
$\|F(\cdot, t)\|_{C^{2+\sigma}(\Sigma_1,\oM)}$
is uniformly bounded for $t\in [0,T)$. Using the
uniformly boundedness of $\|\bar{\nabla}^k B\|_{\bar{g}_+}$ we conclude
uniformly boundedness of
$\|F(\cdot, t)\|_{C^{k+1+\sigma}(\Sigma_1,\oM)}$.
\end{proposition}
\begin{corollary} \label{corollary 6.1}
 $T=+\infty$ and there exists a sequence $t_n\rightarrow +\infty$
such that
$\sup_{\Sigma_1}\|H_{t_n}\|$ $ \to 0,$
 when
$n\to +\infty$.
Moreover, if $\sup_{p\in\Sigma_1}d_2(f_t(p), f_0(p))<c_8$ uniformly
for $t\in[0,+\infty)$, where $c_8>0$ is a constant, then
$f_{t_n}\rightarrow f_{\infty}$ when $n\rightarrow +\infty$,
where $f_{\infty}:\Sigma_1\rightarrow \Sigma_2$ is a smooth totally
geodesic map defining a graphic spacelike
totally geodesic submanifold. Furthermore, if $Ricci_1(p)>0$ at some
point $p\in \Sigma_1$, then
$f_{\infty}$ is constant.
\end{corollary}
\noindent
{\em Proof.}
Note that $F_t(p)$ is a curve in $\oM$ with derivative $H_t$.
Let  $\bar{d}$ be the distance function on $(\oM,\bar{g}_+)$.
We have a uniform bound $\|H_t\|_{g_+}\leq c_9$,
for a constant $c_9>0$, and from (1.1), for $t\geq s$
\begin{equation} \label{eq 6.1}
F_t(p)=F_s(p)+\int_{s}^tH_{\tau}(p)d\tau \mbox{~~~~and~~~~}
\bar{d}(F_t(p),F_s(p))\leq \int_s^t\|H_{\tau}\|_{g_+}d\tau.
\end{equation}
Assuming $T<+\infty$ we have
$\|F_t\|_{L^{\infty}(\Sigma_1,\mathbb{R}^{N})}$  uniformly
bound for $t\in [0,T)$. Therefore,
 by proposition 6.1,  $\|F(\cdot, t)\|_{C^{2+\sigma}(\Sigma_1,
\oM)}$
are uniformly bounded. We take a sequence $t_N\to T$.
By the Ascoli-Arzela theorem we may extract a subsequence $t_n\to T$
of $t_N$, such that
$F(\cdot, t_n)$  converges uniformly to a map $F(\cdot, T)$
in $C^{2}(\Sigma_1,\oM)$, with
$t_n\rightarrow T$.
This also implies  $F(\cdot,t)$ converges uniformly
to $F(\cdot, T)$ when $t\rightarrow T$.
To see this we  only have to note that
for
$d(F_t,F_T)=sup_{p\in \Sigma_1}\bar{d}(F_t(p), F_T(p))$
or $d$ defined from the 2-H\"{o}lder norm,
the following inequality holds:
\begin{equation} \label{eq 6.2}
d(F_t,F_T)\leq d(F_t,F_{t_n})
+d(F_{t_n},F_T)\leq c_{10}|t-t_n|+ d(F_{t_n},F_T),
\end{equation}
where $c_{10}>0$ is a
constant, and in the last inequality we used proposition 5.3.
We note
that $F_T$ is  smooth, by using, by induction,  higher order Schauder
theory to sequential subsequences of $F_{t_n}$, and finally
a diagonal one.
 Following the same reasoning as in the proof of proposition
5.2, $F_{T}$ is a spacelike graph of a map $f_{T}\in
C^{\infty}(\Sigma_1,\Sigma_2)$, and consequently we can extend the solution
$F_t$ to $[0, T+\epsilon)$ for some $\epsilon>0$,
 what is impossible. This proves
$T=+\infty$. It follows from proposition 5.2 that
$$\int_0^{+\infty}\sup_{\Sigma_1}\|H_t\|^2\leq c_{12},$$
for some constant $ c_{12}>0$. Consequently there exist
$t_N\rightarrow +\infty$ with
$\sup_{\Sigma_1}\|H_{t_N}\|^2\rightarrow 0$. Assuming
$f_t$ lies in a compact set of $\Sigma_2$ we are assuming
$\|F_t\|_{L^{\infty}(\Sigma_1,\oM)}\leq C$
 uniformly for $t\in [0, +\infty)$, what implies,
 as above in this proof, for  $t_n$ subsequence of $t_N$,
 $F_{t_n}$ converges
to a map $F_{\infty}\in C^{\infty}(\Sigma_1,\oM)$ when
$n\rightarrow +\infty$, with  $F_{\infty}$  a
spacelike graph of a map $f_{\infty}\in C^{\infty}(\Sigma_1,\Sigma_2)$.
 Now $F_{\infty}$  must satisfy
$\sup_{\Sigma_1}\|H_{\infty}\|^2= 0$,
 that is  $\Gamma_{f_{\infty}}$ is maximal.
From the Bernstein results in \cite{alal,ls} we conclude
$\Gamma_{f_{\infty}}$ is a totally geodesic submanifold, or equivalently,
$f:(\Sigma_1, g_1) \rightarrow (\Sigma_2, g_2)$ is a totally geodesic
map,
and if $Ricci_1>0$ somewhere, then $f_{\infty}$ is a constant map.
 \qed
\\[4mm]
So we have proved parts (1) and  (2) of theorem 1.1.

\section{The case $Ricci_1>0$ everywhere}
\renewcommand{\thesection}{\arabic{section}}
\renewcommand{\theequation}{\thesection.\arabic{equation}}
\setcounter{equation}{0}
\setcounter{theorem}{0}
\setcounter{corollary}{0}
\setcounter{proposition}{0}
\setcounter{lemma}{0}

Next we assume $Ricci_1>0$ everywhere, and prove the last part (3)
of theorem \ref{theorem 1.1},  giving a
particular version of the proof of theorem  \ref{theorem 1.1}, with
no need of using
the Bernstein theorems obtained in \cite{alal,  ls}, but with
a direct proof of  convergence at infinity of
all flow to a graph of a constant map.
\begin{lemma} \label{lemma 7.1}
If $Ricci_1>0$ everywhere, then for $\eta$ given in (\ref{eq 5.3}),
\begin{eqnarray}
1\leq \eta_t^2 \leq 1+ c_{16}e^{-2c_{15}t}, \label{eq 7.1}\\
\lambda_i^2(p,t)\leq \frac{c_{16}e^{-2c_{15} t}}{(1+ c_{16}e^{-2c_{15} t})}
~~~~~~\forall i. \label{eq 7.2}
\end{eqnarray}
for some  constants $c_{15},c_{16}>0$.
 Thus
$\eta_t\rightarrow 1$ when
$t\rightarrow +\infty$.
\end{lemma}
\noindent
{\em Proof.}
The assumption on the sectional
 curvatures of $\Sigma _1$  and $\Sigma_2$ in Theorem 1.1
with the further assumption $Ricci_1>0$ everywhere guarantees that
for each $i$ fixed,
$$\left(\frac{1}{(1-\lambda _i^2)}Ricci_1(a_i,a_i)
+\sum_{j\neq i}\frac{\lambda_j^2}{(1-\lambda_i^2)(1-\lambda_j^2)}
\left[K_1(P_{ij})-K_2(P'_{ij})\right]\right)\geq c_{14}$$
 for some constant $c_{14}>0$. Then
we have by Proposition 4.3 and the proof of the first part of
Theorem 1.1 (see (\ref{eq 5.4}))
$$\frac d{dt}\ln (\cosh\theta)\le \Delta \ln (\cosh\theta)-c_{15}\sum_i
\lambda_i^2,$$ for a constant $c_{15}>0$.
Now recall that $\lambda_i\leq 1$ and
$$1\ge\prod_i(1-\lambda _i^2)=1-\sum _i\lambda _i^2+\sum _{i<j}\lambda
_i^2\lambda _j^2-\sum _{i<j<k}\lambda _i^2\lambda _j^2\lambda
_k^2+\cdots.$$
By induction on $m$ we see that
$$A_m:=\prod_{1\leq i\leq m}(1-\lambda _i^2)- 1+\sum_{1\leq i \leq  m}
\lambda_i^2 = A_{m-1}+\lambda_{m}^2\La{(}1-\prod_{1\leq i\leq m-1}
(1-\lambda _i^2) \La{)} \geq 0.$$
We therefore have,
${\cosh^{-2}\theta}=\prod_i(1-\lambda _i^2)\ge 1-\sum _i\lambda _i^2.$
It follows that for the positive constant $c_{15}$
\begin{eqnarray} \label{eq 7.3}
\frac d{dt}\ln (\cosh\theta)\le \Delta \ln (\cosh\theta)
+c_{15}\left(\frac
1{\cosh^2\theta}-1\right).\end{eqnarray}
 If there exists a time
$t_0$ such that $\eta =1$, then $M_{t_0}$ is already a stable
solution. So without loss of generality, we may assume $\eta \ne
1$.  In view of (\ref{eq 7.3})
and lemma 5.2 we have
$$\frac{d \eta}{dt}\leq c_{15}\left(\frac{1-\eta^2}{\eta}\right),$$
and so
$
1\leq \eta^2\leq 1+c_{16}e^{-2c_{15}t}$
where $c_{16}=\sup_{\Sigma_1}\cosh^2\theta_0 -1$. Consequently,
$$(1-\lambda_i^2)\geq \prod_i(1-\lambda_i^2)
\geq \frac{1}{(1+c_{16}e^{-2c_{15}t})}.~~~ $$\qed

\begin{proposition} \label{proposition 7.1}
If  $Ricci_1>0$ everywhere, then $f_t$ lies
in a compact region of $\Sigma_1$ and
$f_t$ $C^{\infty}$-converges to a unique constant map when
$t\rightarrow +\infty$.
\end{proposition}
\noindent
{\em Proof.}
Set $\varepsilon(t)=c_{16}e^{-2c_{15}t}$ and $\delta(t)=1-\varepsilon(t)$.
By (\ref{eq 7.2}), $\lambda_i^2\leq \varepsilon(t)$.
On the other hand, from (\ref{eq 5.4})
$$||B||^2-\sum
_{i, k}\lambda _i^2(h_{ik}^{m+i})^2-2\sum _{k,i<j}\lambda _i\lambda
_jh_{ik}^{m+j}h_{jk}^{m+i}\ge \delta
(t)||B||^2.
$$
 It follows from Proposition 4.3 that
\begin{eqnarray*}
\frac d{dt}\cosh\theta &\le &\Delta \cosh\theta-\delta
(t)\cosh\theta ||B||^2.
\end{eqnarray*}
Set $\bar\eta =\cosh\theta$. We have by (3.9) $||\nabla \bar\eta
||^2\le \varepsilon (t) \bar\eta ^2||B||^2.$ Let $p(t)$ be a
positive function not less than 1 of $t$, and $\phi=\bar\eta
^p||B||^2$. We have
\begin{eqnarray*}
\frac d{dt}\bar\eta ^p &=&\bar\eta ^p (\frac {dp}{dt}\ln \bar\eta
+p\frac
1{\bar\eta}\frac d {dt}\bar\eta)\\
&=& \bar\eta ^p\frac{dp}{dt}\ln \bar \eta +p\bar\eta ^{p-1}\frac d {dt}\bar\eta \\
&\le & \bar\eta ^p\frac{dp}{dt}\ln \bar \eta+p\bar\eta ^{p-1}(\Delta \bar\eta -\delta\bar\eta ||B||^2)\\
&=&\bar\eta ^p\frac{dp}{dt}\ln \bar \eta+\Delta \bar\eta
^p-p(p-1)\bar\eta ^{p-2}|\nabla \bar\eta |^2
-p\delta\bar\eta ^p||B||^2\\
&\le&\bar\eta ^p\frac{dp}{dt}\ln \bar \eta+\Delta \bar\eta ^p-p
\delta\bar\eta ^p||B||^2,
\end{eqnarray*}
and using (\ref{eq 5.6})
\begin{eqnarray*}
\frac d{dt}\phi &\le& (\bar\eta ^p\frac{dp}{dt}\ln \bar \eta+\Delta
\bar\eta ^p- p\delta\bar\eta ^p||B||^2)||B||^2 \nonumber\\
&&+\bar\eta
^p(\Delta ||B||^2+c_1||B||+c_2||B||^2-\frac 2n||B||^4)\nonumber\\
&=& \Delta \phi-2\bar\eta ^{-p}\nabla \bar\eta ^p\nabla \phi
+2\bar\eta ^{-p}|\nabla \bar\eta ^p|^2||B||^2\\
&&+\bar\eta ^p\left[c_1||B||+(\frac{dp}{dt}\ln \bar
\eta+c_2)||B||^2-(p\delta +\frac 2n)||B||^4\right].
\end{eqnarray*}
Now   $2\bar\eta ^{-p}\|\nabla \bar\eta ^p\|^2||B||^2\le
2p^2\varepsilon \bar\eta ^p||B||^4,
$ and we have
\begin{eqnarray*}
\frac d{dt}\phi &\le& \Delta \phi-2\bar\eta ^{-p}\nabla \bar\eta
^p\nabla \phi\\
&& +\bar\eta ^p\left[c_1||B||+(\frac{dp}{dt}\ln \bar
\eta+c_2)||B||^2+(2p^2\varepsilon -p\delta-\frac
2n)||B||^4\right]\nonumber\\
&=&\Delta \phi-2\bar\eta ^{-p}\nabla \bar\eta ^p\nabla \phi
+\bar\eta ^{-p}(2p^2\varepsilon -p\delta-\frac 2n)\phi^2\\
&&+c_1\bar\eta
^{\frac p2}\phi^{\frac 12}+(\frac{dp}{dt}\ln \bar \eta+c_2)\phi.
\end{eqnarray*}
Setting $\psi =e^{\frac 12 c_{15}t}\phi $, then $\psi$ satisfies the
evolution inequality
\begin{eqnarray*}
\frac d{dt}\psi &\le& \Delta \psi-2\bar\eta ^{-p}\nabla \bar\eta
^p\nabla \psi +e^{-\frac 12c_{15}t}\bar\eta ^{-p}(2p^2\varepsilon
-p\delta-\frac 2n)\psi^2\\
&&+c_1e^{\frac 14c_{15}t}\bar\eta ^{\frac
p2}\psi^{\frac 12}+(\frac{dp}{dt}\ln \bar \eta+c_2+\frac 12c_{15})\psi.
\end{eqnarray*}
Taking $p^2=\frac 1{n\varepsilon}$,  that is,
$$p(t)=\frac{e^{c_{15}t}}{
\sqrt{nc_6}},$$
and using that $\ln(1+a)\leq a$ for all $a\geq 0$, we obtain
$\ln \bar\eta\le \frac {c_{16}}{2}e^{-2c_{15}t}.$
It is easy to check that
$$(1+c_{16}e^{-2c_{15}t})
^{\frac p4}=(1+c_{16}e^{-2c_{15}t}) ^{\frac
{e^{c_{15}t}}{4\sqrt{nc_{16}}}}\rightarrow 1, \mbox{ as }t\rightarrow
\infty,$$
what implies $1\geq \bar{\eta}^{-p}\geq c_{13}$, for some constant
$c_{13}>0$.
Therefore
\begin{eqnarray*}
\frac d{dt}\psi &\le& \Delta \psi-2\bar\eta ^{-p}\nabla \bar\eta
^p\nabla \psi
 -\frac{c_{13}\delta}{\sqrt{nc_{16}}}e^{\frac 12 c_{15}t}\psi^2\\
&&+c_1e^{\frac 14c_{15}t}(1+c_{16}e^{-2c_{15}t}) ^{\frac p4}\psi^{\frac
12}+(\frac{c_{15}\sqrt{c_{16}}}{2\sqrt{n}}e^{-c_{15}t}+c_2+\frac 12c_{15})\psi.
\end{eqnarray*}
and by choosing $T_0$ large enough  such that for $t\ge T_0$ one has
$\delta(t)\ge \frac 12$ and
\begin{eqnarray} \label{eq 7.4}
\frac d{dt}\psi &\le& \Delta \psi-2\bar\eta ^{-p}\nabla \bar\eta
^p\nabla \psi -c_{17}\left\{e^{\frac 12 c_{15}t}\psi^2-e^{\frac
14c_{15}t}\psi^{\frac 12}-\psi\right\},
\end{eqnarray}
for a constant $c_{17}>0$.\\[2mm]
{\bf Claim:} When $t\ge T_0$, $||B||^2\le c_{19}e^{-\tau t}$ for some
positive constants $c_{19}$ and $\tau$.
\\[2mm]
We prove the claim.
 For any $t_0\in [T_0, +\infty)$, we consider
for each $t\in [T_0, t_0]$ a point $x_t$ such that
 $\psi (x_t,t)$ attains its maximum $\psi_{\max}(t)$.
Since $\psi_{\max}$ is a locally Lipschitz function on $[T_0, t_0]$,
we  may take $t_1\in [T_0, t_0]$ a point where this maximum is
achieved. If $t_1=T_0$, we have done. Thus we may assume $t_1>T_0$.
At $(x_1,t_1)$, \ $\Delta \psi\leq 0$, $\nabla \psi=0$ and $\frac
d{dt}\psi\geq 0$. Thus from (\ref{eq 7.4}) at $(x_{t_1},t_1)$,
$$e^{\frac 12 c_{15}t}\psi^2-e^{\frac
14c_{15}t}\psi^{\frac 12}-\psi\le 0.$$
Let $\bar\psi (x_{t_1},t_1)=\sqrt{\psi (x_{t_1}, t_1)}$. Then,
$e^{\frac{1}{2}c_{15}t_1}\bar{\psi}^3-\bar{\psi}\leq e^{\frac{1}{4}c_{15}t_1}$,
what implies at $(x_{t_1},t_1)$
$$\bar{\psi}^3\leq e^{\frac{1}{4}c_{15}t_1}\bar{\psi}^3
\leq 1+e^{-\frac{1}{4}c_{15}t_1}\bar{\psi}\leq 1+\bar{\psi}.$$ Thus,
there is a constant $c_{18}>0$ that does not depend on $t_0$, such
that $\bar\psi(x_{t_1},t_1)\le c_{18}$. Therefore,
$$\max_{M_{t_0}}\psi = \psi(x_{t_0},t_{0})\leq
\psi(x_{t_1},t_{1}) \le c_{18}^2.$$
 Then for a constant $c_{19}$ not depending on $t_0$ and $t$,
$\max _{M_{t_0}}||B||^2\le c_{19}e^{-\frac
12 c_{15}t_0}.$ Since $t_0$ is arbitrary, we prove the claim.\\[2mm]
By the claim, we have the estimate of $||H||$
\begin{equation} \label{eq 7.5}
||H||^2\le \frac {c_{19}}{m}e^{-\tau t}.
\end{equation}
By (\ref{eq 7.5}) and (\ref{eq 6.1}),
we have $\bar{d}(F_t(p), F_0(p))\le c_{20}$ for
a positive constant
$c_{20}$. Therefore $F_t(p)$ lies in a compact region of $\bar M$.
From the proof of corollary \ref{corollary 6.1},
$F_{t_n}$ converges uniformly to a limit
map $F_{\infty}$ for a sequence
 $t_n\rightarrow \infty$. Moreover $F_{\infty}$ is
a spacelike graph of a map $f_{\infty}$. By (\ref{eq 7.2})
$f_{\infty}=constant$.
Now we prove the limit is unique. Following the proof of corollary
\ref{corollary 6.1}, in (\ref{eq 6.2}) we have in this case
$${d}(F_t,F_{\infty})\leq  \left|\int_{t}^{t_n}  \|H\|_{\max}ds\right|
+ d(F_{t_n},F_{\infty})
\leq  \frac{c_{19}\tau}{m}|e^{-\tau t}- e^{-\tau t_n}| +d(F_{t_n},F_{\infty})$$
and the right hand side converges to zero.
This implies $F_t$ to converge to $F_{\infty}$
in $C^{\infty}(\Sigma_1, \oM)$. To see this, note that if
there exists $t_N\to +\infty$ such that $d(F_{t_N},F_{\infty})\geq c_{30}$,
$c_{30}$ a positive constant, where $d$ is the distance function
relative to the $C^1$- H\"{o}lder norm,
 then, extracting a subsequence
$t_n$ of $t_N$ with $F_{t_n}$ converging in $C^1(\Sigma_1,\oM)$
this sequence must converge to $F_{\infty}$, leading to a contradiction.
By induction we prove $F_t$ converges to $F_{\infty}$ in $C^k$,
for each $k\geq 0$.
 \qed\\[4mm]
\noindent
 {\em  Proof of theorem 1.2}.
We consider the
pseudo-Riemannian product space $\oM=\Sigma_1\times \Sigma_2$
equipped with the pseudo-Riemannian metric
$g_1-g'_2$, $g'_2=\rho^{-1} g_2$ for a constant $\rho >0$. 
The  spacelike condition for $f$ means $f^*g_2<\rho g_1$.
The curvature tensor $R'$ of $g'_2$ satisfies
$R'_2(X',Y', X',Y')=\rho R_2(X,Y,X,Y)$, where $X'=\sqrt{\rho}X$
and $Y'=\sqrt{\rho}Y$. If $K_1>0$, then
$K_1\geq {K'}_2$ is satisfied if we assume $\rho \leq 
{\min_{\Sigma_1} K_1}/{\sup_{\Sigma_2}K_2^+}$. 
If $K_2\leq 0$  we may take $\rho=+\infty$.
If $Ricci_1>0$ and $K_2\leq -c<0$, we may take any
$\rho\geq {\max_{\Sigma_1}K_1^-}/c$, where $K^-=\max\{-K,0\}$,
and in particular $\rho=+\infty$.
Then we are in the conditions of theorem 1.1(3).
\qed \\[4mm]
{\em Proof of corollary 1.2. }
Under the assumptions, $\sum_{i}\lambda_i^2+ 1\leq
\Pi_i(1+\lambda_i^2)<2$. In particular, $\Gamma_f$ is a spacelike graph for
the corresponding pseudo-Riemannian structure of $\Sigma_1\times
\Sigma_2$. \qed\\[4mm]
{\em Proof of corollary 1.3. }
We assume $f$ minimizes the $\phi$-energy functional.
Let $f_t$ given by  theorem 1.2. From
(7.2) of lemma 7.1 and the assumptions on $\phi$,
we have   $E_{\phi}(f)\leq \liminf E_{\phi}(f_t)= E_{\phi}(f_{\infty})=0$, when
$t\to +\infty$. Thus, $E_{\phi}(f)=0$, and so $f$ is constant.
\qed\\[4mm]
{\em Remark. }
We note that when $E_{\phi}$ is the usual energy functional of $f$  
corollary 1.3 can be obtained 
using a simple  Weitzenb\"{o}ck formula.
Since $f$ is harmonic
$$\Delta\|df\|^2=\|\nabla df\|^2+\sum_{i\neq j} 
\lambda_i^2 \left( K_1(P_{ij})-\lambda_j^2K_2(P'_{ij})\right)$$
that under the curvature conditions of theorem 1.2,
$K_1\geq K_2$ with $K_1\geq 0$, or $Ricci_1\geq 0$
and $K_2\leq 0$, we 
have $\Delta\|df\|^2\geq 0$, what implies $\|df\|$, and so by the above
equation $f$ is totally geodesic, and if $Ricci_1>0$ at some point,
then $\lambda_i=0$, that is $f$ is constant.
So, the corollary is mainly interesting for $\phi$ not the square
of the norm. But this argument also shows that the curvature condition
$Ricci_1\geq 0$ with $K_2<0$  is
the expected one, since in this case, 
 it is well known \cite{es}, that
any map $f:\Sigma_1\to \Sigma_2$ is homotopic to an harmonic map,
and so, under the condition $Ricci_1> 0$,
 necessarily to a constant one.

 \end{document}